\documentclass{article}
\usepackage{latexsym,amssymb,amsmath}
\usepackage{graphicx}
\usepackage{color}
\usepackage[emtex,matrix,arrow]{xy}

\newtheorem{thm}{Theorem}[section]
\newtheorem{cor}[thm]{Corollary}
\newtheorem{ex}[thm]{Example}
\newtheorem{lem}[thm]{Lemma}
\newtheorem{prop}[thm]{Proposition}
\newtheorem{defn}[thm]{Definition}
\newtheorem{rem}[thm]{Remark}
\numberwithin{equation}{section}
\newcommand{\Real}{\mathbb R}
\newcommand{\Natural}{\mathbb N}

\newcommand{\To}{\longrightarrow}

\newcommand{\Aa}{\mathcal{A}}

\newcommand{\Cc}{\mathcal{C}}
\newcommand{\Bb}{\mathcal{B}}
\newcommand{\Dd}{\mathcal{D}}
\newcommand{\F}{\mathcal{F}}
\newcommand{\G}{\mathcal{G}}

\newcommand{\Cgg}{\mathfrak{C}}
\newcommand{\Dgg}{\mathfrak{D}}
\newcommand{\id}{\rm{id}}
\font\Gotica= eufm10

\def\Cg{{\Gotica\char 67}}
\def\Dg{{\Gotica\char 68}}

\newcommand{\qed}{\hspace*{\fill}$\Box$  \ifmmode \else
    \par\addvspace\topsep\fi}
\newenvironment {proof}{\par\addvspace\topsep\noindent{\it Proof.}
    \ignorespaces }{\qed}

\begin{document}

\title{2-cosemisimplicial objects in a 2-category, permutohedra and deformations
of pseudofunctors}
\author{Josep Elgueta \\
Dept. Matem\`atica Aplicada II \\ Universitat Polit\`ecnica de
  Catalunya \\ email: Josep.Elgueta@upc.es}



\maketitle

\begin{abstract}
In this paper we take up again the deformation theory for
$K$-linear pseudofunctors initiated in \cite{jE1}. We start by
introducing a notion of 2-cosemisimplicial object in an arbitrary
2-category and analyzing the corresponding coherence question,
where the permutohedra make their appearence. We then describe a
general method to obtain usual cochain complexes of $K$-modules
from (enhanced) 2-cosemisimplicial objects in the 2-category ${\bf
Cat}_K$ of small $K$-linear categories and prove that the
deformation complex $X^{\bullet}(\F)$ introduced in \cite{jE1} can
be obtained by this method from a 2-cosemisimplicial object that can
be associated to $\F$. Finally, using this
2-cosemisimplicial object of $\F$ and a generalization to the context of
$K$-linear categories of the deviation calculus introduced by Markl
and Stasheff for $K$-modules \cite{MS94}, it is shown that the
obstructions to the integrability of an $n^{th}$-order deformation
of $\F$ indeed correspond to cocycles in the third cohomology group
$H^3(X^{\bullet}(\F))$, a question which remained open in
\cite{jE1}.
\end{abstract}


\section{Introduction}

In \cite{jE1}, we introduced a deformation complex for $K$-linear
unitary pseudofunctors which turned out to describe the so-called
purely pseudofunctorial first order deformations. This was a
generalization to the many objects setting of Yetter's deformation
theory for monoidal functors (see \cite{dY98},\cite{dY01}). A
common feature of both deformation theories, which also appears in
other categorical or 2-categorical deformation theories, such as
Crane and Yetter's deformation theory for semigroupal categories
\cite{CY981}) or the deformation theory for semigroupal
2-categories \cite{jE1}, is the presence of suitable ``padding
operators'' in the definition of the coboundary maps. These
operators may look like something artificial in the construction.
One of the purposes of this paper is to give a framework where they
appear most naturally. Our point of view is that the presence of
such padding operators is a consequence of the intrinsically
higher-dimensional nature of the structures that are being
deformed. Conjecturally, they are the shadow of a
higher-dimensional description, still to be found, of the
corresponding deformation theory. In this sense, we guess that the
right setting for studying categorical deformations should involve
a suitable notion of 2-cochain complex, together with the
corresponding notion of 2-co(semi)simplicial object in a 2-category.
Along these lines, we introduce in this paper a notion of
2-cosemisimplicial object in an arbitrary 2-category (a
2-dimensional version of the classical cosemisimplicial objects in
a category), and we show that the deformation complex of a
$K$-linear unitary pseudofunctor $\F$ can be obtained from such an
object that may be associated to $\F$. It is precisely in this
process of going from the 2-cosemisimplicial object to the cochain
complex that the padding operators appear. Presumably, this process
involves a loss of information. It is then tempting to think that
more information should be contained in the hypothetical {\it
2-cochain complex} that should be derived from the
2-cosemisimplicial object, and that this 2-cochain complex could
give a more complete description of the deformations of the
pseudofunctor (including, for example, deformations at the level of
1-morphisms). At this point, it is worth mentioning the works by R.
Street on cohomology with coefficients in an ($n$-)category
\cite{rS87}, \cite{rS95}, \cite{rS03}. This author has
recently given (see \cite{rS03}) a precise definition of what he
calls the {\it
  descent $n$-category} of any cosimplicial $n$-category
$\mathcal{E}^{\bullet}$. It seems possible that this notion of
descent $n$-categories (or some variant of it) provides the right setting we
are claiming for to give the cohomological description of the deformations of
higher dimensional algebraic structures.

As in any categorification process, in defining the notion of
2-cosemisimplicial object in a 2-category, suitable coherence
conditions are introduced and the corresponding coherence theorem
should be proved. In doing this, it turns out that the
{\it permutohedra}, first introduced by Milgram in the context of
iterated loop spaces \cite{rjM66}, are the right family of convex
polytops describing the higher-order cosemisimplicial
identities, in a way analogous to that encountered when weakening
the associativity equation, where the role is played by the famous
Stasheff associahedra.

The last purpose of the paper concerns higher-order obstructions.
It remained as an open question in \cite{jE1} if the obstructions
to the integrability of an $n^{th}$-order deformation indeed live
in one of the cohomology groups, a condition which, according to
Gerstenhaber \cite{mG64}, must satisfy any good cohomological
deformation theory. We prove that this is indeed the case. More
explicitly, we show that the obstructions correspond to 3-cocycles
in the deformation complex introduced in \cite{jE1}. To prove this,
we use a generalization to the context of $K$-linear categories of
the Markl and Stasheff deviation calculus \cite{MS94}. As it will
be seen, the previously constructed 2-cosemisimplicial object turns
out to be quite useful in making the proof easy to write.

The paper is organized as follows. Section 2 contains some
definitions and preliminary results needed later. In Section 3, we
recall the notion of deformation of a pseudofunctor we work with as
well as the definition of the deformation complex as given in
\cite{jE1}. In Section 4 we define 2-cosemisimplicial objects in an
arbitrary (strict) 2-category and prove the corresponding coherence
theorem. In Section 5 we focus the attention on the special case
of the 2-category ${\bf Cat}_K$ of (small) $K$-linear categories
and show that in this case usual cochain complexes of $K$-modules
can be obtained from a suitably {\it enhanced} 2-cosemisimplicial
object in ${\bf Cat}_K$. In the next section, we go back to the
deformation theory of a pseudofunctor, proving that one can
construct a (trivially enhanced) 2-cosemisimplicial object from any
pseudofunctor and that, when the pseudofunctor is $K$-linear, its
deformation complex coincides with one of the cochain complexes one
may obtain by the method in the previous section. Finally, in
Section 7 we generalize Markl and Stasheff deviation calculus to
the context of arbitrary $K$-linear categories. This technique is
used in the next section to prove that the obstructions to the
integrability of a partial deformation indeed live in the
corresponding cohomology.


\section{Preliminaries}

Unless otherwise indicated, $K$ denotes a given commutative field. Let us first
recall the definition of a pseudofunctor
between 2-categories (see, for ex., \cite{fB94}).

%
%
%
%
%
%
%

\begin{defn} \label{pseudofunctor}
If \Cg\ and \Dg\ are two 2-categories, a {\sl pseudofunctor} from
$\Cgg$ to $\Dgg$ is any quadruple $\F=(|\F|,\F_*,\widehat{\F}_*,\F_0)$, where
\begin{itemize}
\item
$|\F|:|\Cgg|\rightarrow|\Dgg|$ is an object map;

\item
$\F_*=\{\F_{X,Y}:\Cgg(X,Y)\rightarrow\Dgg(\F(X),\F(Y))\}$
is a collection of functors, indexed by ordered pairs of objects
$X,Y\in|\Cgg|$;

\item
$\widehat{\F}_*=\{\widehat{\F}_{X,Y,Z}:
c^{\Dgg}_{\F(X),\F(Y),\F(Z)}\circ(\F_{X,Y}\times\F_{Y,Z})
\Rightarrow\F_{X,Z}\circ c^{\Cgg}_{X,Y,Z}\}$
is a collection of natural isomorphisms,
indexed by ordered triples of objects $X,Y,Z\in|\Cgg|$ (the $c^{\Cgg}_{-,-,-}$
denote the composition functors in the 2-category $\Cgg$ and similarly
$c^{\Dgg}_{-,-,-}$). Explicitly,
this means having a 2-isomorphism \footnote{In this paper, the arguments in
  $\widehat{\F}$
  are written in the reverse order to that used in \cite{jE1}.}
$$
\widehat{\F}_{X,Y,Z}(f,g):\F_{Y,Z}(g)\circ\F_{X,Y}(f)\Rightarrow\F_{X,Z}(g\circ
f)
$$
for any path of 1-morphisms
$X\stackrel{f}{\rightarrow} Y\stackrel{g}{\rightarrow} Z$,
natural in $(f,g)$, and

\item
$\F_0=\{\F_0(X):\F_{X,X}(\id_X)\Rightarrow \id_{\F(X)}\}$ is a
collection of 2-isomorphisms, indexed by objects $X\in|\Cgg|$.
\end{itemize}
\noindent{Moreover}, this data must satisfy the following {\sl
  coherence axioms} (for short, the indexing objects are omitted so that we
just write $\widehat{\F}(f,g)$ and $\F(f)$):

\begin{itemize}
\item[(A1)]
({\sl Composition axiom}) For all paths of 1-morphisms
$X\stackrel{f}{\rightarrow} Y\stackrel{g}{\rightarrow}
Z\stackrel{h}{\rightarrow} T$, the following diagram commutes
$$
\xymatrix{
\F(h)\circ\F(g)\circ\F(f)\ar@{=>}[d]_{\widehat{\F}(g,h)\circ
  1_{\F(f)}}
\ar@{=>}[rr]^{1_{\F(h)}\circ\widehat{\F}(f,g)} &
& \F(h)\circ\F(g\circ f)\ar@{=>}[d]^{\widehat{\F}(h,g\circ f)} \\
\F(h\circ g)\circ\F(f)\ar@{=>}[rr]_{\widehat{\F}(f,h\circ g)} & &
\F(h\circ g\circ f) }
$$

\item[(A2)]
({\sl Unit axioms}) For any 1-morphism $f:X\rightarrow
Y$, the following equalities hold:
\begin{align*}
\widehat{\F}({\id}_X,f)&=1_{\F(f)}\circ\F_0(X) \\
\widehat{\F}({\id}_Y,f)&=\F_0(Y)\circ 1_{\F(f)}
\end{align*}
\end{itemize}
\end{defn}

The whole set of 2-isomorphisms $\widehat{\F}(f,g)$ and $\F_0(X)$, for
all objects $X$ and
composable 1-morphisms $f,g$, will be called the {\sl pseudofunctorial
  structure} on $\F$. When they are all identities
the pseudofunctor is called a {\sl 2-functor}. When only the
$\F_0(X)$ are identities, we will call it a {\sl unitary
pseudofunctor}.

For later use, we give in the next Lemma a ``component-free''
description of the above composition axiom. The proof is
an easy exercise left to the reader.

\begin{lem} \label{axioma_composicio_reescrit}
Let $\F=(|\F|,\F_*,\widehat{\F}_*,\F_0)$ be the data defining a pseudofunctor
between two 2-categories $\Cgg$ and $\Dgg$, and let us define families of
functors \footnote{The meaning of the notation used to distinguish these
  families will be seen in Section 6.}
$\{\F^{(1,1,1)}_{X,Y,Z,T}\}$,
$\{\F^{(1,2)}_{X,Y,Z,T}\}$, $\{\F^{(2,1)}_{X,Y,Z,T}\}$
and $\{\F^{(3)}_{X,Y,Z,T}\}$ and natural isomorphisms
  $\{\sigma^{12}_{X,Y,Z,T}\}$, $\{\sigma^{24}_{X,Y,Z,T}\}$,
  $\{\sigma^{13}_{X,Y,Z,T}\}$ and $\{\sigma^{34}_{X,Y,Z,T}\}$, both indexed by
  ordered quadruples $(X,Y,Z,T)$ of objects
  in $\Cgg$, and respectively given by
\begin{align}
&\begin{array}{l}
\F^{(1,1,1)}_{X,Y,Z,T}=c^{\Dgg}_{\F(X),\F(Z),\F(T)}\circ(c^{\Dgg}_{\F(X),\F(Y),\F(Z)}\times
\id_{\Dgg(\F(Z),\F(T))})\circ \\ \hspace{5 truecm}
\circ(\F_{X,Y}\times\F_{Y,Z}\times\F_{Z,T}) \end{array}
\label{functor_associat_1}
\\
&\F^{(1,2)}_{X,Y,Z,T}=c^{\Dgg}_{\F(X),\F(Y),\F(T)}\circ(\F_{X,Y}\times\F_{Y,T})
\circ(\id_{\Cgg(X,Y)}\times
c^{\Cgg}_{Y,Z,T}) \\ &\F^{(2,1)}_{X,Y,Z,T}=c^{\Dgg}_{\F(X),\F(Z),\F(T)}\circ(\F_{X,Z}\times\F_{Z,T})\circ(c^{\Cgg}_{X,Y,Z}\times
\id_{\Cgg(Z,T)}) \\ &\F^{(3)}_{X,Y,Z,T}=\F_{X,T}\circ
c^{\Cgg}_{X,Z,T}\circ(c^{\Cgg}_{X,Y,Z}\times{\id}_{\Cgg(Z,T)}) \label{functor_associat_4}
\end{align}
and \footnote{In this paper, identity 2-morphisms are generically denoted by
  $1_f$. But when the 1-morphism $f$ is a functor we use a boldface 1, to
  emphasize the fact that it is an identity natural transformation.}
\begin{align}
\sigma^{12}_{X,Y,Z,T}&=\mathbf{1}_{c^{\Dgg}_{\F(X),\F(Z),\F(T)}}\circ
(\widehat{\F}_{X,Y,Z}\times\mathbf{1}_{\F_{Z,T}}) \label{trans_associada_1}\\
\sigma^{24}_{X,Y,Z,T}&=\widehat{\F}_{X,Z,T}\circ\mathbf{1}_{c^{\Cgg}_{X,Y,Z}\times
id_{\Cgg(Z,T)}} \\ \sigma^{13}_{X,Y,Z,T}&=\mathbf{1}_{c^{\Dgg}_{\F(X),\F(Y),\F(T)}}\circ
(\mathbf{1}_{\F_{X,Y}}\times\widehat{\F}_{Y,Z,T})
\\ \sigma^{34}_{X,Y,Z,T}&=\widehat{\F}_{X,Z,T}\circ\mathbf{1}_{id_{\Cgg(X,Y)}\times
c^{\Cgg}_{Y,Z,T}} \label{trans_associada_4}
\end{align}
Then, the previous composition axiom is equivalent to the commutativity of the
diagrams of natural transformations
\begin{equation} \label{axioma_hexagon_inicial}
\xymatrix{
\F^{(1,1,1)}_{X,Y,Z,T}\ar@{=>}[rr]^{\sigma^{12}_{X,Y,Z,T}}
\ar@{=>}[d]_{\sigma^{13}_{X,Y,Z,T}} & &
\F^{(2,1)}_{X,Y,Z,T}\ar@{=>}[d]^{\sigma^{24}_{X,Y,Z,T}} \\
\F^{(1,2)}_{X,Y,Z,T}\ar@{=>}[rr]_{\sigma^{34}_{X,Y,Z,T}} & &
\F^{(3)}_{X,Y,Z,T} }
\end{equation}
for all ordered quadruples $(X,Y,Z,T)$ of objects in $\Cgg$.
\end{lem}

The above definitions may be extended to the $K$-linear context using
the Deligne product between $K$-linear categories and functors (see, for ex.,
\cite{dY01}, Chap. 10). Furthermore,
we will need
to define the $K[[h]]$-linear extensions of the corresponding $K$-linear
versions. Such definitions already appear in \cite{jE1}, although they were
formulated without using the notion of Deligne product.

Recall that by a $K$-{\it linear category} one
means a category $\Cc$ enriched over the monoidal category
\boldmath $Vect_K$ \unboldmath of $K$-vector spaces. The corresponding
topological version will be
called a {\it complete $K[[h]]$-linear category}. By definition, it is a
category
\boldmath$\Cc$\unboldmath\ enriched over the monoidal category
\boldmath$K[[h]]$-$Mod_c$\unboldmath\ of separated and complete
$K[[h]]$-modules.

\begin{defn}
A $K$-{\it linear 2-category} is a 2-category
\Cg\ whose hom-categories $\Cgg(X,Y)$, for all objects $X,Y$ of $\Cgg$, are
$K$-linear, and whose composition functors
$c^{\Cgg}_{X,Y,Z}:\Cgg(X,Y)\times\Cgg(Y,Z)\To\Cgg(X,Z)$, for all $X,Y,Z$,
are $K$-bilinear or, equivalently, $K$-linear functors
$c^{\Cgg}_{X,Y,Z}:\Cgg(X,Y)\odot\Cgg(Y,Z)\To\Cgg(X,Z)$, where $\odot$
denotes the Deligne product of $K$-linear categories.

Similarly, a {\it complete $K[[h]]$-linear 2-category} is a
2-category \boldmath $\Cgg$ \unboldmath whose hom-categories $\mbox{\boldmath$\Cgg$}(X,Y)$, for all objects $X,Y$ of
\boldmath $\Cgg$ \unboldmath, are
complete $K[[h]]$-linear categories and whose composition functors
$c^{\Cgg}_{X,Y,Z}:\mbox{\boldmath$\Cgg$}(X,Y)\times\mbox{\boldmath$\Cgg$}(Y,Z)\To
\mbox{\boldmath$\Cgg$}(X,Z)$,
for all $X,Y,Z$, are $K[[h]]$-bilinear or, equivalently,
$K[[h]]$-linear functors
$c^{\Cgg}_{X,Y,Z}:\mbox{\boldmath$\Cgg$}(X,Y)\widehat{\odot}\mbox{\boldmath$\Cgg$}(Y,Z)
\To\mbox{\boldmath$\Cgg$}(X,Z)$,
where $\widehat{\odot}$ denotes the topological Deligne product
of complete $K[[h]]$-linear categories.
\end{defn}

\begin{defn}
Given two $K$-linear 2-categories $\Cgg,\Dgg$, a $K$-linear
pseudofunctor from $\Cgg$ to $\Dgg$ is a pseudofunctor
$\F:\Cgg\To\Dgg$ whose defining functors
$\F_{X,Y}:\Cgg(X,Y)\To\Dgg(F(X),F(Y))$, for all objects $X,Y$ of
$\Cgg$, are $K$-linear.

Similarly, by replacing the term $K$-linear by (complete)
$K[[h]]$-linear, one gets the definition of $K[[h]]$-linear pseudofunctor
between complete $K[[h]]$-linear 2-categories.
\end{defn}

Notice that the defining natural isomorphisms
$\widehat{\F}_{X,Y,Z}:c^{\Dgg}_{\F(X),\F(Y),\F(Z)}\circ(\F_{X,Y}\times\F_{Y,Z})\Rightarrow\F_{X,Z}\circ
c^{\Cgg}_{X,Y,Z}$ of a $K$-linear pseudofunctor $\F$
may also be considered as natural transformations
$\widehat{\F}_{X,Y,Z}:c^{\Dgg}_{\F(X),\F(Y),\F(Z)}\circ(\F_{X,Y}\odot\F_{Y,Z})\Rightarrow
\F_{X,Z}\circ c^{\Cgg}_{X,Y,Z}$. The same thing is true for a $K[[h]]$-linear
pseudofunctor,
with the topological Deligne product $\widehat{\odot}$ replaced
by $\odot$. The reader may easily check that there is also an analog of
Lemma~\ref{axioma_composicio_reescrit}
for $K$-linear and complete $K[[h]]$-linear pseudofunctors, where the
cartesian product $\times$ in the definition of
the
functors (\ref{functor_associat_1})-(\ref{functor_associat_4}) and natural
transformations (\ref{trans_associada_1})-(\ref{trans_associada_4}) must be
replaced by the
Deligne product $\odot$ and the topological Deligne product $\widehat{\odot}$,
respectively.

We will be mostly concerned with the $K[[h]]$-linear extensions of
a $K$-linear 2-category or pseudofunctor. Let us first recall the
definitions in the context of categories.

Given a $K$-linear category $\Cc$, its $K[[h]]$-linear
extension, denoted by $\Cc[[h]]$, is the complete
$K[[h]]$-linear category with the same objects as $\Cc$ and
$K[[h]]$-modules of morphisms given by
$$
\Cc_h(X,Y):=\Cc(X,Y)[[h]]\ , \quad X,Y\in|\Cc|
$$
where $A[[h]]$, for any $K$-module $A$, denotes the topologically free
$K[[h]]$-module of
formal power series in $h$ with coefficients in $A$.
Composition in $\Cc[[h]]$ is defined
in the obvious way in terms of the composition in $\Cc$ and the
product rule of formal power series. In particular, the identity morphisms in
$\Cc[[h]]$ are the same as in $\Cc$.
It seems that these categories were introduced for the first time by Drinfeld
\cite{vD90} in his study of the quasiHopf algebras, providing the setting for
the deformation theory of monoidal categories (see Crane and Yetter
\cite{CY981} and Yetter \cite{dY01}). For its later use, let us state the
following result, whose
proof is left to the reader (it
is the analog in the context of categories of a well-known result about the
topological tensor product between topologically free $K[[h]]$-modules):

\begin{lem} \label{iso_producteDeligne_categories}
For any $K$-linear categories $\Cc$, $\Dd$ there is an isomorphism of complete
$K[[h]]$-linear categories
$\Psi_{\Cc,\Dd}:\Cc[[h]]\widehat{\odot}\Dd[[h]]\stackrel{\cong}{\To}
(\Cc{\odot}\Dd)[[h]]$.
\end{lem}

Given a $K$-linear functor $F:\Cc\To\Dd$ between
$K$-linear categories, its $K[[h]]$-linear extension, denoted by
$F[[h]]$, is the $K[[h]]$-linear functor
$F[[h]]:\Cc[[h]]\To\Dd[[h]]$ acting on objects as $F$ and such that
$$
F[[h]]\left(\sum_{k\geq 0}f_kh^k\right)=\sum_{k\geq 0}F(f_k)h^k
$$
It is easy to check that $(F'\circ F)[[h]]=F'[[h]]\circ F[[h]]$ and
$(id_{\Cc})[[h]]=id_{\Cc[[h]]}$ for all composable $K$-linear functors $F,F'$
and $K$-linear categories $\Cc$. The proof of the next lemma is also left to
the reader.

\begin{lem} \label{lema_igualtat}
For any $K$-linear categories $\Cc_1,\Cc_2,\Dd_1,\Dd_2$ and $K$-linear
functors $F_i:\Cc_i\To\Dd_i$, $i=1,2$, we have
$$
(F_1\odot F_2)[[h]]\circ\Psi_{\Cc_1,\Cc_2}=\Psi_{\Dd_1,\Dd_2}\circ
(F_1[[h]]\widehat{\odot}F_2[[h]])
$$
\end{lem}
Another easy but important fact needed later is the following:

\begin{lem} \label{lema_nat_h}
For any $K$-linear functors $F,G:\Cc\To\Dd$ between arbitrary
$K$-linear categories $\Cc,\Dd$, there is an isomorphism of $K[[h]]$-modules
$$
{\rm Nat}(F,G)[[h]]\cong{\rm Nat}(F[[h]],G[[h]])
$$
sending the formal power series $\sum_{k\geq 0}\tau_kh^k$ to the natural
transformation $\tau_h:F[[h]]\Rightarrow G[[h]]$ with components
$$
\left(\sum_{k\geq 0}\tau_kh^k\right)_X=\sum_{k\geq 0}(\tau_k)_Xh^k,\qquad
X\in|\Cc|
$$
Furthermore, under this identification, the vertical and horizontal
compositions of naturals transformations are given by the usual product rule
of formal power series.
\end{lem}

\begin{proof}
By definition, a natural transformation $\tau_h:F[[h]]\Rightarrow G[[h]]$
involves a collection of morphisms $(\tau_h)_X:F(X)\To G(X)$ in
$\Dd[[h]]$, for all objects $X$ of $\Cc$. But a generic such morphism
is of the form
$$
(\tau_h)_X=\sum_{n\geq 0}(\tau_n)_Xh^n
$$
The proof reduces to show that the naturality of $(\tau_h)_X$
in $X$ is equivalent to the naturality in $X$ of the $(\tau_n)_X$,
for all $n\geq 0$. This last condition may be shown by an easy
induction which is left to the reader. As regards the formula for
the vertical composition, it immediately follows from the
definition of composition in $\Dd[[h]]$.
\end{proof}

The corresponding notions of $K[[h]]$-linear extension in the 2-category
setting can now be formulated as follows.

\begin{defn}
Let \Cg\ be a $K$-linear 2-category. Then, its $K[[h]]$-{\it linear
extension} is the complete $K[[h]]$-linear 2-category $\Cgg[[h]]$
given by the following data:
\begin{itemize}
\item[(i)]
The objects of $\Cgg[[h]]$ are the same as in $\Cgg$.
\item[(ii)]
The hom-categories $\Cgg[[h]](X,Y)$ are the $K[[h]]$-linear extensions of the
corresponding categories, i.e., for all objects $X,Y$,
\begin{equation} \label{def_Cg_h}
\Cgg[[h]](X,Y):=\Cgg(X,Y)[[h]].
\end{equation}
\item[(iii)]
The composition functors
$c^{\Cgg[[h]]}_{X,Y,Z}:\Cgg[[h]](X,Y)\widehat{\odot}\Cgg[[h]](Y,Z)\To\Cgg[[h]](X,Z)$,
for any objects $X,Y,Z$ of $\Cgg$, are given by
\begin{equation} \label{def_Cg_h_b}
c^{\Cgg[[h]]}_{X,Y,Z}:=c^{\Cgg}_{X,Y,Z}[[h]]\circ\Psi_{X,Y,Z}
\end{equation}
where $c^{\Cgg}_{X,Y,Z}[[h]]$ is the $K[[h]]$-linear
extension of the composition functor $c^{\Cgg}_{X,Y,Z}$ of $\Cgg$ and
$\Psi_{X,Y,Z}=\Psi_{\Cgg(X,Y),\Cgg(Y,Z)}$ (see
Lemma~\ref{iso_producteDeligne_categories}).
\item[(iv)]
The identity 1-morphisms $id_X$ are the same as in $\Cgg$.
\end{itemize}
\end{defn}

The reader may easily check that the above data indeed defines a
(K[[h]]-linear) 2-category. Notice that according to
(\ref{def_Cg_h}), the 1-morphisms in $\Cgg[[h]]$ are exactly the
same as in $\Cgg$ but a generic 2-morphism $\tau_h:f\Rightarrow f'$
between two such 1-morphisms $f,f':X\To Y$ is of the form of a
formal power series
$$
\tau_h=\tau_0+\tau_1 h+\tau_2 h^2+\cdots
$$
with the $\tau_i:f\Rightarrow f'$, $i\geq 0$, 2-morphisms in \Cg.
Also implicit in (\ref{def_Cg_h}) is the fact that the vertical
composition of two such 2-morphisms is given by the usual product rule
of formal power series, while
(\ref{def_Cg_h_b}) means that the composition of
1-morphisms in $\Cgg[[h]]$ is the same as in $\Cgg$ and the
horizontal composition of two 2-morphisms $\tau_h:f\Rightarrow
f':X\To Y$ and $\sigma_h:g\Rightarrow g':Y\To Z$ is given by the
product rule.

Before giving the corresponding notion of $K[[h]]$-linear extension
for $K$-linear pseudofunctors, let us first remark that for any
$K$-linear pseudofunctor $\F:\Cgg\To\Dgg$, we have (see Lemma
\ref{lema_igualtat})
\begin{align*}
c^{\Dgg[[h]]}_{\F(X),\F(Y),\F(Z)}&\circ(\F_{X,Y}[[h]]
\widehat{\odot}\F_{Y,Z}[[h]])= \\ =(&c^{\Dgg}_{\F(X),\F(Y),\F(Z)}\circ(\F_{X,Y}
\odot\F_{Y,Z}))[[h]]\circ \Psi_{X,Y,Z}
\end{align*}
and
$$
\F_{X,Z}[[h]]\circ c^{\Cgg[[h]]}_{X,Y,Z}=(\F_{X,Z}\circ
c^{\Cgg}_{X,Y,Z})[[h]]\circ\Psi_{X,Y,Z}
$$
Hence, the following definition makes sense (see also Lemma~\ref{lema_nat_h}).

\begin{defn} \label{extensio_pseudof}
Let $\F:\Cgg\To\Dgg$ be a $K$-linear pseudofunctor between
$K$-linear 2-categories. Then, the $K[[h]]$-linear extension of
$\F$ is the $K[[h]]$-linear pseudofunctor
$\F[[h]]:\Cgg[[h]]\To\Dgg[[h]]$ acting on objects as $\F$ and whose
remaining structural data is given by:
\begin{itemize}
\item[(i)]
$\F[[h]]_{X,Y}=\F_{X,Y}[[h]]$ (the $K[[h]]$-linear extension of $\F_{X,Y}$).
\item[(ii)]
$\widehat{\F[[h]]}_{X,Y,Z}=\widehat{\F}_{X,Y,Z}\circ
{\mathbf 1}_{\Psi_{X,Y,Z}}$
(here, $\widehat{\F}_{X,Y,Z}$ stands for a formal power
series of natural transformation with only zero order
term).
\item[(iii)]
$\F[[h]]_0(X)=\F_0(X)$.
\end{itemize}
\end{defn}

We leave to the reader to check that the previous data indeed
define a $K[[h]]$-linear pseudofunctor between $\Cgg[[h]]$ and
$\Dgg[[h]]$. Notice that, according to conditions (i) and (ii), for
any path $X\stackrel{f}{\To}Y\stackrel{g}{\To}Z$ and any 2-morphism
$\tau_h=\tau_0+\tau_1h+\cdots:f\Rightarrow f'$ in
$\Cgg[[h]]$ we have
\begin{align*}
\F[[h]](f)&=\F(f) \\
\F[[h]](\tau_h)&=\F(\tau_0)+\F(\tau_1)h+\cdots \\
\widehat{\F[[h]]}(f,g)&=\widehat{\F}(f,g)
\end{align*}


\section{Deformation complex of a $K$-linear pseudofunctor}

Given a $K$-linear pseudofunctor $\F$, we introduced in \cite{jE1} a cochain
complex $(X^{\bullet}(\F),\delta)$ which in the unitary case described
the purely pseudofunctorial first order deformations of $\F$. A fundamental question
which remained open was if the
obstructions to the integrability of a partial deformation live in some of the
cohomology groups. This point is settled down in
Section 8 using an analog of Markl and Stasheff deviation calculus \cite{MS94}. In this section, we recall the
necessary definitions from \cite{jE1}.


\begin{defn}
Let $\Cgg,\Dgg$ two $K$-linear 2-categories and $\F:\Cgg\To\Dgg$ a
$K$-linear pseudofunctor. Then, by a purely pseudofunctorial formal
deformation of $\F$ we mean any $K[[h]]$-linear pseudofunctor
$\F_h:\Cgg[[h]]\To\Dgg[[h]]$ differing from the $K[[h]]$-linear
extension $\F[[h]]$ (see Definition~\ref{extensio_pseudof}) only in the
pseudofunctorial structure, which must be of the form
\begin{equation} \label{F_barret}
(\widehat{\F}_h)_{X,Y,Z}=\left(\sum_{k\geq
0}\widehat{\F}^k_{X,Y,Z}h^k\right)\circ {\mathbf 1}_{\Psi}
\end{equation}
\begin{equation} \label{F_0}
(\F_h)_0(X)=\sum_{k\geq 0}\F^k_0(X)h^k
\end{equation}
with $\widehat{\F}^0_{X,Y,Z}=\widehat{\F}_{X,Y,Z}$ and
$\F^0_0(X)=\F_0(X)$ for all objects $X,Y,Z$ of \Cg.
\end{defn}
Notice that $\F[[h]]$ itself gives an example of such a deformation,
called the {\it null deformation}, where $\widehat{\F}^k_{X,Y,Z}=0$ and
$\F^k_0(X)=0$ for all $k\geq 1$.

Clearly, a purely pseudofunctorial formal
deformation of $\F$ is completely given by the families
of natural transformations $\{\widehat{\F}^k_{X,Y,Z}\}_{X,Y,Z}$ and
2-morphisms $\{\F^k_0(X)\}_X$, for all $k\geq 1$.
However, they are not arbitrary. They must
be such that the corresponding natural transformations
(\ref{F_barret}) and 2-morphisms (\ref{F_0}) indeed define a pseudofunctorial
structure on $\F_h$. Next result makes precise the conditions they
must satisfy in a form suitable to our purposes. In particular, the diagrams
which appear are of the
right kind for the notion of deviation introduced in Section 7 to make
sense.

\begin{lem} \label{condicions_deformacio}
Let $\F:\Cgg\To\Dgg$ be a $K$-linear pseudofunctor. Then, the families
$\{\widehat{\F}^k_{X,Y,Z}\}_{X,Y,Z}$ and $\{\F^k_0(X)\}_X$, $k\geq 1$, define
a purely pseudofunctorial formal deformation $\F_h$ of $\F$ if and only if
\begin{itemize}
\item[(1)]
For all objects $X,Y,Z,T\in|\Cgg|$, the following diagram commutes:
\begin{equation} \label{axioma_hexagon}
\xymatrix{
\F^{(1,1,1)}_{X,Y,Z,T}[[h]]\ar@{=>}[rr]^{\sigma^{12}_{X,Y,Z,T}(h)}
\ar@{=>}[d]_{\sigma^{13}_{X,Y,Z,T}(h)} & &
\F^{(2,1)}_{X,Y,Z,T}[[h]]\ar@{=>}[d]^{\sigma^{24}_{X,Y,Z,T}(h)} \\
\F^{(1,2)}_{X,Y,Z,T}[[h]]\ar@{=>}[rr]_{\sigma^{34}_{X,Y,Z,T}(h)} & &
\F^{(3)}_{X,Y,Z,T}[[h]] }
\end{equation}
where
\begin{align}
\sigma^{12}_{X,Y,Z,T}(h)&:=\sum_{k\geq
0}\left[\mathbf{1}_{c^{\Dgg}_{\F(X),\F(Z),\F(T)}}\circ
(\widehat{\F}^k_{X,Y,Z}\odot\mathbf{1}_{\F_{Z,T}})\right]h^k \label{tau_1} \\
\sigma^{24}_{X,Y,Z,T}(h)&:=\sum_{k\geq
0}\left[\widehat{\F}^k_{X,Z,T}\circ\mathbf{1}_{c^{\Cgg}_{X,Y,Z}\odot
\id_{\Cgg(Z,T)}}\right]h^k \\ \sigma^{13}_{X,Y,Z,T}(h)&:=\sum_{k\geq
0}\left[\mathbf{1}_{c^{\Dgg}_{\F(X),\F(Y),\F(T)}}\circ
(\mathbf{1}_{\F_{X,Y}}\odot\widehat{\F}^k_{Y,Z,T})\right]h^k
\\ \sigma^{34}_{X,Y,Z,T}(h)&:=\sum_{k\geq
0}\left[\widehat{\F}^k_{X,Z,T}\circ\mathbf{1}_{\id_{\Cgg(X,Y)}\odot
c^{\Cgg}_{Y,Z,T}}\right]h^k \label{tau_2}
\end{align}

\item[(2)]
For all objects $X,Y\in|\Cgg|$, all 1-morphisms $f:X\To Y$ and all $k\geq 1$,
the following equalities hold
\begin{eqnarray}
&\widehat{\F}^k_{X,X,Y}({\id}_X,f)=1_{\F(f)}\circ\F_0^k(X) \label{triangle1} \\
&\widehat{\F}^k_{X,Y,Y}(f,{\id}_Y)=\F_0^k(Y)\circ 1_{\F(f)} \label{triangle2}\end{eqnarray}
\end{itemize}
\end{lem}
The set of equations (\ref{axioma_hexagon}) together with
(\ref{triangle1})-(\ref{triangle2}) play the role of the associativity
equation in the study of the formal deformations of an associative algebra
\cite{mG64}, and are called the {\it structural} or
{\it deformation equations}.
\begin{proof}
By the topological $K[[h]]$-linear version of
Lemma~\ref{axioma_composicio_reescrit}, we know that the composition axiom is
equivalent to the commutativity of the diagrams
\begin{equation} \label{diagrama_desviacio=obstruccio}
\xymatrix{
(\F_h)^{(1,1,1)}_{X,Y,Z,T}\ar@{=>}[rr]^{(\sigma_h)^{12}_{X,Y,Z,T}}
\ar@{=>}[d]_{(\sigma_h)^{13}_{X,Y,Z,T}} & &
(\F_h)^{(2,1)}_{X,Y,Z,T}\ar@{=>}[d]^{(\sigma_h)^{24}_{X,Y,Z,T}} \\
(\F_h)^{(1,2)}_{X,Y,Z,T}\ar@{=>}[rr]_{(\sigma_h)^{34}_{X,Y,Z,T}} & &
(\F_h)^{(3)}_{X,Y,Z,T} }
\end{equation}
for all objects $X,Y,Z,T\in|\Cgg|$. Now,
using Lemma~\ref{lema_igualtat}, we obtain that
%
$$
(\F_h)^{\alpha}_{X,Y,Z,T}=\F^{(\alpha)}_{X,Y,Z,T}[[h]]\circ \Psi_{X,Y,Z,T}
$$
for all $\alpha=(1,1,1),(2,1),(1,2),(3)$, where
$\Psi_{X,Y,Z,T}=\Psi_{\Cgg(X,Y),\Cgg(Y,Z),\Cgg(Z,T)}$ denotes the canonical
isomorphism
$\Cgg(X,Y)[[h]]\widehat{\odot}\Cgg(Y,Z)[[h]]\widehat{\odot}\Cgg(Z,T)[[h]]\cong(\Cgg(X,Y)\odot\Cgg(Y,Z)\odot\Cgg(Z,T))[[h]]$,
whose existence follows from Lemma~\ref{iso_producteDeligne_categories}.
On the other hand, a straightforward computation shows that
$$
(\sigma_h)^{ij}_{X,Y,Z,T}=\sigma^{ij}_{X,Y,Z,T}(h)\circ{\mathbf
    1}_{\Psi_{X,Y,Z,T}}
$$
for all pairs $i,j$, where the $\sigma^{ij}_{X,Y,Z,T}(h)$ are the natural
transformations
(\ref{tau_1})-(\ref{tau_2}). Hence, condition
(\ref{axioma_hexagon_inicial}) on $\F_h$ takes the form
$$
\xymatrix{
\F^{(1,1,1)}_{X,Y,Z,T}[[h]]\circ\Psi_{X,Y,Z,T}
\ar@{=>}[rrr]^{\sigma^{12}_{X,Y,Z,T}(h)\circ{\mathbf 1}_{\Psi_{X,Y,Z,T}}}
\ar@{=>}[d]_{\sigma^{13}_{X,Y,Z,T}(h)\circ{\mathbf 1}_{\Psi_{X,Y,Z,T}}} & & &
\F^{(2,1)}_{X,Y,Z,T}[[h]]\circ\Psi_{X,Y,Z,T}
\ar@{=>}[d]^{\sigma^{24}_{X,Y,Z,T}(h)\circ{\mathbf 1}_{\Psi_{X,Y,Z,T}}} \\
\F^{(1,2)}_{X,Y,Z,T}[[h]]\circ\Psi_{X,Y,Z,T}
\ar@{=>}[rrr]_{\sigma^{34}_{X,Y,Z,T}(h)\circ{\mathbf 1}_{\Psi_{X,Y,Z,T}}} & & &
\F^{(3)}_{X,Y,Z,T}[[h]]\circ\Psi_{X,Y,Z,T} }
$$
By the interchange law this is equivalent to
$$
(\sigma^{24}_{X,Y,Z,T}(h)\cdot\sigma^{12}_{X,Y,Z,T}(h))\circ{\mathbf 1}_{\Psi_{X,Y,Z,T}}=
(\sigma^{34}_{X,Y,Z,T}(h)\cdot\sigma^{13}_{X,Y,Z,T}(h))\circ{\mathbf 1}_{\Psi_{X,Y,Z,T}}
$$
and, since $\Psi_{X,Y,Z,T}$ is an isomorphism (in particular, essentially
surjective),
the terms in $\Psi_{X,Y,Z,T}$ may indeed be cancelled to give the equivalent
condition (\ref{axioma_hexagon}). The proof that equalities
(\ref{triangle1})-(\ref{triangle2}) are in turn equivalent to the unit
axioms on the deformed pseudofunctor $\F_h$ is left to
the reader.
\end{proof}
Together with the notion of purely pseudofunctorial {\it formal} deformation,
in \cite{jE1} we also introduced the corresponding notion of purely
pseudofunctorial {\it $n^{th}$-order} deformation, for all $n\geq 1$. It
is defined in the same way as the formal deformations by replacing the ring
of formal power series $K[[h]]$ by the ring of
truncated
polynomials $K[h]/(h^n)$. Using arguments similar to those made above, it may
be shown that such a deformation is
completely given by families $\{\widehat{\F}^k_{X,Y,Z}\}$ and $\{\F^k_0(X)\}$
as
above, for $k=1,\ldots,n$,  satisfying the deformation equations
(\ref{axioma_hexagon}) up to $h^{n+1}$ and (\ref{triangle1})-(\ref{triangle2})
for all $k=1,\ldots,n$. The details are left to the reader.

Then, for any $K$-linear pseudofunctor
$\F:\Cgg\rightarrow\Dgg$, we defined in \cite{jE1} a cochain complex
$X^{\bullet}(\F)$
whose vector spaces $X^n(\F)$ were given by
$$
X^n(\F):=\left\{
\begin{array}{ll} \prod_{(X_0,\ldots,X_n)\in|\Cgg|^{n+1} }
{\rm
  Nat}(\F^{(1,\stackrel{n)}{\ldots},1)}_{X_0,\ldots,X_n},
\F^{(n)}_{X_0,\ldots,X_n}) & n\geq 1 \\
0 & {\rm otherwise}
\end{array}\right.
$$
where
\begin{align*}
\F^{(1,\stackrel{n)}{\ldots},1)}_{X_0,\ldots,X_n}&:=c^{\Dgg}_{\F(X_0),\ldots,\F(X_n)}
\circ(\F_{X_0,X_1}\odot\F_{X_1,X_2}\odot\cdots\odot\F_{X_{n-1},X_n})
\\
\F^{(n)}_{X_0,\ldots,X_n}&:=\F_{X_0,X_n}\circ c^{\Cgg}_{X_0,\ldots,X_n}
\end{align*}
for all $n\geq 2$ (they are the components of two particular $\F$-iterates of
multiplicity $n$ chosen as references) and
$$
\F^{(1)}_{X_0,X_1}:=\F_{X_0,X_1}
$$
if $n=1$ (the unique $\F$-iterate of multiplicity 1). Here, the $c^{\Cgg}$ and
$c^{\Dgg}$ indexed by $n+1$  objects, $n\geq 3$, denote the unique
$n^{th}$-order induced composition functors in the corresponding
2-category. The coboundary map $\delta:X^{n-1}(\F)\longrightarrow
X^n(\F)$, $n\geq 2$, was then defined in terms of the ``padding''
operators $\lceil-\rceil$ associated to $\F$ (see \cite{jE1}) by the formula
\begin{align*}
(\delta\phi)(f_0,f_1,\ldots,f_{n-1})=&\lceil
1_{\F(f_{n-1})}\circ\phi(f_0,\ldots,f_{n-2})
\rceil_{\F(X_0),\F(X_n)}
\\ &+\sum_{i=1}^{n-1}(-1)^i\lceil\phi(f_0,\ldots,f_i\circ f_{i-1},\ldots,f_{n-1})\rceil_{\F(X_0),\F(X_n)}
\\ &+(-1)^n\lceil\phi(f_1,\ldots,f_{n-1})\circ
1_{\F(f_0)}\rceil_{\F(X_0),\F(X_n)}
\end{align*}
with $\phi\in X^{n-1}(\F)$ and
$f_i\in|\Cgg(X_i,X_{i+1})|$, $i=0,\ldots,n-1$ (notice that 1-morphisms $f_i$
are indexed differently with respect to the notation in
\cite{jE1} and that, as
arguments of $\phi$, they are written in the reverse order).
We proved then the following:

\begin{thm} \cite{jE1}
Let $\F$ be a $K$-linear unitary pseudofunctor and let us denote by
$H^{\bullet}(\F)$ the
cohomology of the corresponding deformation complex as defined above. Then,
the equivalence classes of the purely pseudofunctorial first order
deformations are in one-one
correspondence with the elements of $H^2(X^{\bullet}(\F))$.
\end{thm}


\section{2-cosemisimplicial objects in a 2-category}

As mentioned in the introduction, in this section we introduce a
notion of 2-cosemisimplicial object in a 2-category as a sort of
categorification of the classical notion of cosemisimplicial object
in a category (see, for ex., \cite{cW94}). Our original motivation for
doing this was to see that,
associated to any pseudofunctor between 2-categories, we have such
an object, and that the cochain complex of a $K$-linear
pseudofunctor in the previous section can be obtained from it. This
is done in Section 6.

Recall that, given any category $\Cc$, a cosemisimplicial object in
$\Cc$ is any covariant functor $K:\Delta_s\To\Cc$, where $\Delta_s$
(the semisimplicial category) is the subcategory of the simplicial
category $\Delta$ whose morphisms are the injections
$\alpha:[i]\hookrightarrow [n]$ (see \cite{cW94}). To define the
corresponding categorified notion, $\Cc$ should be replaced by a
bicategory $\Cgg$, $\Delta_s$ by a suitable `semisimplicial
  bicategory' $2\Delta_s$ and $K:\Delta_s\To\Cc$ by a pseudofunctor
$\F:2\Delta_s\To\Cgg$. The outstanding point is what we should take
as semisimplicial bicategory $2\Delta_s$. A priori, the only
reasonable condition we have on it is that it should be a categorification of
$\Delta_s$. But the categorification of a given mathematical
structure is not unique in general. For example, the set $\Natural$ of
natural numbers as a ``rig'' has the category of finite sets as well
as the category of finite 
dimensional vector spaces over a given field $K$ as two nonequivalent
categorifications, or the usual
notion of commutative monoid, which has both the notions of symmetric
monoidal category and braided monoidal 
category as two nonequivalent categorifications. To avoid making such
a choice and at 
the same time to have a description as explicit as possible of the
notion of 2-cosemisimplicial object, we will take as our starting
point the definition of cosemisimplicial object in $\Cc$ which
follows from the presentation of $\Delta_s$ in terms of generators
and relations. Thus, using such presentation of $\Delta_s$, it can
be shown that a cosemisimplicial object in $\Cc$ is the same thing
as a sequence of objects $K_0,K_1,\cdots$ in $\Cc$ together with
{\it coface morphisms} $\partial^i_n:K^{n-1}\To K^n$,
$i=0,\ldots,n$, $n\geq 1$, satisfying the {\it cosemisimplicial
identities}
\begin{equation} \label{identitats_cosemisimplicials}
\partial^j_{n+1}\circ\partial^i_n=\partial^i_{n+1}\circ\partial^{j-1}_n,
\qquad 0\leq i<j\leq n+1
\end{equation}
We then take as definition in the 2-dimensional setting the
following (to simplify, we further restrict to the context of
2-categories).

\begin{defn} \label{objecte_2-cosemisimplicial}
Given a 2-category $\Cgg$, a 2-cosemisimplicial object in $\Cgg$ is
any sequence of objects $X^0,X^1,\ldots$ in $\Cgg$ together with
1-morphisms (the {\sl coface maps}) $\partial^i_n:X^{n-1}\To X^n$,
for all $i=0,\ldots,n$ and $n\geq 1$, and 2-isomorphisms (the {\sl
 cosemisimplicial coherers}) $\tau_{ij}^n:\partial^j_{n+1}\circ\partial^i_n\Rightarrow\partial^i_{n+1}\circ\partial^{j-1}_n$,
$0\leq i<j\leq n+1$, such that the diagrams
$$
\xymatrix{
\partial^k_{n+2}\circ\partial^j_{n+1}\circ\partial^i_n\ar@{=>}[rr]^{1_{\partial^k_{n+2}}\circ\\
  \tau_{ij}^n}
\ar@{=>}[d]^{\tau_{jk}^{n+1}\circ 1_{\partial^i_n}} & &
\partial^k_{n+2}\circ
\partial^i_{n+1}\circ\partial_n^{j-1}\ar@{=>}[rr]^{\tau^{n+1}_{ik}\circ
1_{\partial^{j-1}_n}} & &
\partial^i_{n+2}\circ\partial^{k-1}_{n+1}\circ\partial^{j-1}_n\ar@{=>}[d]_
{1_{\partial^i_{n+2}}\circ\tau^n_{j-1,k-1}}
\\
\partial^j_{n+2}\circ\partial^{k-1}_{n+1}\circ\partial^i_n\ar@{=>}[rr]_
{1_{\partial^j_{n+2}}\circ\tau^n_{i,k-1}} & &
\partial^j_{n+2}\circ\partial^i_{n+1}\circ\partial^{k-2}_n\ar@{=>}[rr]_
{\tau^{n+1}_{ij}\circ 1_{\partial^{k-2}_n}} & & \partial^i_{n+2}\circ\partial^{j-1}_{n+1}\circ\partial^{k-2}_n }
$$
commute for all $0\leq i<j<k\leq n+2$ and all $n\geq 1$.
\end{defn}

For short, such a 2-cosemisimplicial object will be denoted by the
triple $(X^{\bullet},\partial,\tau)$ or just by $X^{\bullet}$, when
there is no confusion. Notice that this definition includes as
special cases the usual cosemisimplicial objects in a category
$\Cc$ when we think of $\Cc$ as the 2-category with only the
identity 2-morphisms.

The commutative diagrams in the above definition are the coherence
laws that appear in any categorification process, and they are
imposed to get the corresponding coherence theorem. To state this
theorem, let us consider, for any $s,k\geq 1$, the subcategory
$\Cc_{s,k}$ of $\Cgg(X^{s-1},X^{s+k})$ with objects all composites
of the coface maps, i.e., all 1-morphisms $f:X^{s-1}\rightarrow
X^{s+k}$ of the form
$$
f=\partial^{i_k}_{s+k}\circ\partial^{i_{k-1}}_{s+k-1}\circ\cdots\circ\partial^{i_0}_s
$$
for $i_j=0,1,\ldots,s+j$ ($j=0,\ldots,k$). We will refer to such
1-morphisms as the $\partial$-{\it paths} from $X^{s-1}$ to
$X^{s+k}$. Given two such $\partial$-paths $f,f'$, the morphisms
from $f$ to $f'$ in $\Cc_{s,k}$ are all possible pastings of the
coherers $\tau^n_{ij}$'s and the identity 2-morphisms of the coface
maps giving a 2-morphism between them. They will be denoted by
$\sigma:f\Rightarrow f'$ because they are actually 2-morphisms in
$\Cgg$. Thus, a generic morphism $\sigma:f\Rightarrow f'$ in
$\Cc_{s,k}$ is of the form
$$
\sigma=(1_{f'_1}\circ\tau^{n_1}_{i_1j_1}\circ 1_{f_1})\cdot
(1_{f'_2}\circ\tau^{n_2}_{i_2j_2}\circ
1_{f_2})\cdot\cdots\cdot(1_{f'_q}\circ\tau^{n_q}_{i_qj_q}\circ
1_{f_q}),
$$
for some $\partial$-paths $f_{\alpha},f'_{\alpha}$ and indices
$i_{\alpha},j_{\alpha},n_{\alpha}$, with $\alpha=1,\ldots,q$ (the dot
denotes the vertical composition of 2-morphisms in $\Cgg$). The
2-morphisms in $\Cgg$ of the form $1_{f'}\circ\tau^n_{ij}\circ 1_f$, for $f,f'$
$\partial$-paths, will be called {\it expanded coherers}. For
example, the
composites $\partial^3_3\circ\partial^1_2\circ\partial^0_1$ and
$\partial^0_3\circ\partial^0_2\circ\partial^1_1$ define two objects of
$\Cc_{1,2}$, and a morphism
in $\Cc_{1,2}$ between them is given by the pasting
$$
(\tau^2_{01}\circ
1_{\partial^1_1})\cdot(1_{\partial^1_3}\circ\tau^1_{02})\cdot(\tau_{13}^2\circ
1_{\partial^0_1})
$$
The coherence theorem states then the following:

\begin{thm} \label{teorema_coherencia}
Let $s,k\geq 1$. Then, for any two objects $f,f'$ in $\Cc_{s,k}$, there is at
most one morphism (actually, an isomorphism) in $\Cc_{s,k}$ from $f$ to $f'$.
\end{thm}
Such a unique isomorphism will be called the {\sl canonical 2-isomorphism}
from $f$ to $f'$, to distinguish it from all other
possible 2-morphisms between $f$ and $f'$ that may exist in $\Cgg$.

To prove the theorem, let us consider the graph $G_{s,k}$ with vertices
all $\partial$-paths $f:X^{s-1}\rightarrow X^{s+k}$ and with edges all the
expanded coherers (hence, $\Cc_{s,k}$ is the
quotient of the free groupoid generated by $G_{s,k}$ modulo the above coherence
relations). It has $(s+1)(s+2)\cdots(s+k+1)$ vertices and it is a degree $k$
regular graph (i.e., for
any vertex, the total
number of incident edges is equal to $k$). It follows that $G_{s,k}$ has
$\frac{1}{2}k(s+1)\cdots(s+k+1)$
edges. Let us identify the vertex
$\partial^{i_k}_{s+k}\circ\partial^{i_{k-1}}_{s+k-1}\circ\cdots\circ\partial^{i_0}_s$
in $G_{s,k}$ with the $(k+1)$-tuple $(i_0,\ldots,i_k)$. The sum
$i_0+\cdots+i_k$ will be called the {\it height} of the
vertex and denoted by $h(i_0,\ldots,i_k)$. We further define the {\it rank} of
such a vertex, denoted by $r(i_0,\ldots,i_k)$, as the number of strictly
positive jumps we meet when going from $i_0$ to
$i_k$. Hence, $0\leq r(i_0,\ldots,i_k)\leq k$. For example,
$r(1,2,3,2,4)=3$ and $r(1,1,2,3)=2$. If we agree that an edge goes {\it out}
of a vertex when the
vertex is the domain of the expanded coherer represented by that edge, while
it goes {\it into}
a vertex when the vertex is its codomain (equivalently, the domain of the
inverse morphism), then the rank of a vertex corresponds to the number of
edges going out of the vertex. A vertex $(i_0,\ldots,i_k)$ will be called
an {\it   out-vertex} when its rank is $k$ (all edges go out of the vertex),
and an {\it in-vertex} when its rank is zero (all edges go into the
vertex). Note that the out-vertices in $G_{s,k}$ are in one-one correspondence
with the subsets of $k+1$ elements of the set $\{0,1,\ldots,s+k\}$, because it
must be $i_0<i_1<\cdots<i_k$. In particular, two differents out-vertices
have different heights. Finally, if the edges of a path in $G_{s,k}$, taken in
order, involve
only expanded coherers and none of its inverses (resp. only inverses of the
expanded coherers), the path will be called {\it
  directed} (resp. {\it inversely directed}).

\begin{figure}[h]
\centering
\input{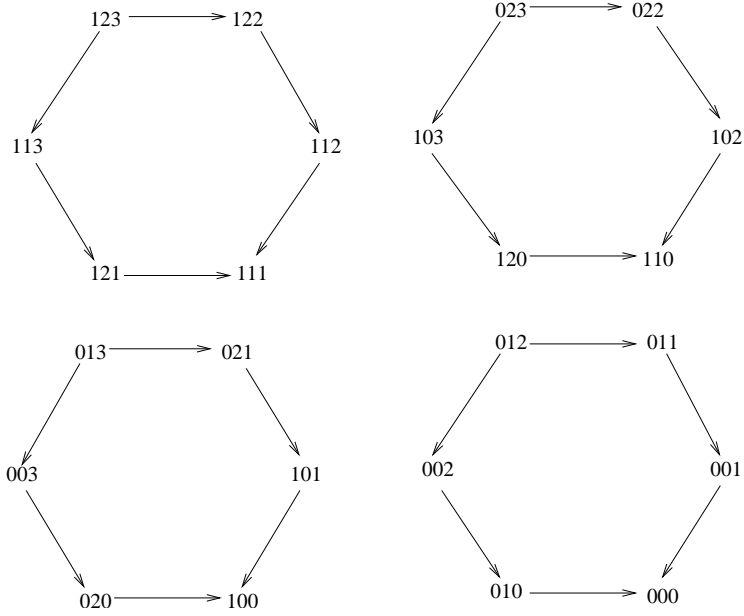}
\caption{The graph $G_{1,2}$}
\label{graf_G12}
\end{figure}

The graph $G_{1,2}$ is depicted in Fig~\ref{graf_G12}. It may be
seen that it has various connected components, all of them isomorphic
and each one with exactly one out-vertex and exactly one
in-vertex. This turns out to be true for all graphs 
$G_{s,k}$, $s,k\geq 1$. To see that, the following
property of $G_{s,k}$ will be used.

\begin{lem} \label{propietat_basica_Gsk}
Let $(i_0,\ldots,i_k)$ be an arbitrary out-vertex in the graph
$G_{s,k}$. Then, any directed path in $G_{s,k}$ from $(i_0,\ldots,i_k)$ to an
in-vertex has length $k(k+1)/2$.
\end{lem}

\begin{proof}
Let us identify each entry $i_p$ ($p=0,\ldots,k$) with its initial
position $p$ in the $(k+1)$-tuple. As we move along a path in
$G_{s,k}$ that starts in this vertex, the $p^{th}$-entry will
change its value (to new values $i'_p$) and the position it
occupies in the $(k+1)$-tuple. The lemma follows from the fact that
we will get an in-vertex when and only when, for any pair of
entries $p<q$, $i'_p$ is to the right of $i'_q$. Indeed, suppose
that, after several edges, there is a pair $p<q$ such that the
$p^{th}$-entry $i'_p$ is still to the left of the $q^{th}$-entry
$i'_q$. We then have $i'_p\in
\{i_p-p,\ldots,i_p\}$, with $i'_p=i_p-p$ when all entries $i_0,\ldots,i_{p-1}$
have been moved to the right of $i_p$, and $i_p'=i_p$ when
none of these entries has been moved to the right of $i_p$. Suppose
$i'_p=i_p-t$ ($t\in\{0,\ldots,p\}$). In this case, we necessarily
have $i'_q\in\{i_q-(q-1-p+t),\ldots,i_q\}$, because $q-1-p+t$ is the
maximum number of positions that $i_q$ can move to the left always
keeping to the right of $i'_p$. It follows that
$$
i'_q-i'_p\geq
i_q-(q-p-1+t)-i_p+t=i_q-i_p-q+p+1\geq q-p-q+p+1=1
$$
and, hence, the vertex is still not an in-vertex. On the other hand, it
is clear that, when all such ``transpositions'' have been made, the
resulting vertex is really an in-vertex. Now, there are $k(k+1)/2$ such
``transpositions'' to be made. Since going through one
directed edge in the graph corresponds to making exactly one of
these ``transpositions'', we conclude that we get an in-vertex
after going over a directed path of length $k(k+1)/2$ and only in
this case.
\end{proof}
Using this lemma, we can prove the following result which will be
used below to prove the coherence theorem, and which in
particular shows that the connected components of $G_{s,k}$ are
parametrized by the injections
$\{0,1,\ldots,k\}\hookrightarrow\{0,1,\ldots,s+k\}$, so that $G_{s,k}$
has $\binom{s+k+1}{k+1}$ connected components. 

\begin{prop} \label{components_connexes_Gnk}
Let $s,k\geq 1$. Then, each connected component of $G_{s,k}$ has
exactly one out-vertex and one in-vertex. Furthermore, all its
components are isomorphic and independent of $s$.
\end{prop}

\begin{proof}
Clearly, each component has at least one out-vertex (just follow
an inversely directed path from any vertex in the component until
the end). To prove that it has at most one, suppose there are two
different out-vertices $(i^{out}_0,\ldots,i^{out}_k)$ and
$(i^{out'}_0,\ldots,i^{out'}_k)$ in the same component $C$. In
particular, they have different heights.  Since there is no
directed path connecting them (no directed path ends in an
out-vertex), there must be directed paths $\gamma$, $\gamma'$
starting at each out-vertex which meet in some common vertex
$(i_0,\ldots,i_k)$. Following from this vertex a directed path
$\overline{\gamma}$ until the end, we will get an in-vertex
$(i^{in}_0,\ldots,i^{in}_k)$. Now, by the previous proposition, all
directed paths from an out- to an in-vertex have the same length,
so that both composite paths $\gamma\overline{\gamma}$ and
$\gamma'\overline{\gamma}$ have the same length. On the other hand,
when going over any directed edge, the height always decreases by
exactly one unit. It follows that the height of the final in-vertex
should have two different values, which makes no sense. Hence, there
is exactly one out-vertex in each component. It immediately follows
then that there is also exactly one in-vertex in
each component, with a well-defined value of its height, equal to
the height of the corresponding out-vertex minus $k(k+1)/2$. To see
that all connected components are isomorphic, let us denote by
$C(i_0^{out},\ldots,i_k^{out})$, $C(i_0^{out'},\ldots,i_k^{out'})$
the connected components corresponding to the out-vertices
$(i_0^{out},\ldots,i_k^{out})$ and
$(i_0^{out'},\ldots,i_k^{out'})$, respectively. Then, for any
vertex $(i_0,\ldots,i_k)$ in $C(i_0^{out},\ldots,i_k^{out})$, we
have
$(i_0,\ldots,i_k)=\mbox{\boldmath$\tau$}(i^{out}_0,\ldots,i^{out}_k)$,
for a suitable composite {\boldmath$\tau$} of expanded coherers.
Then, we get an isomorphism
$\varphi:C(i^{out}_0,\ldots,i^{out}_k)\cong
C(i^{out'}_0,\ldots,i^{out'}_k)$ by defining
$$
\varphi(i^{out}_0,\ldots,i^{out}_k)=(i^{out'}_0,\ldots,i^{out'}_k)
$$
and for any other vertex
$$
\varphi(\mbox{\boldmath$\tau$}(i^{out}_0,\ldots,i^{out}_k))=
\mbox{\boldmath$\tau'$}(i^{out'}_0,\ldots,i^{out'}_k)
$$
where {\boldmath$\tau'$} is the composite of expanded coherers
obtained from {\boldmath$\tau$} by suitably changing the indices of
the expanded coherers which appear in {\boldmath$\tau$}, according
to the corresponding initial out-vertex. Finally, to prove that the
components are independent of $s$, it is enough to see, for ex.,
that the connected components $C(1,\ldots,k+1)$ of $G_{s,k}$ and
$G_{s',k}$, for any $s,s'\geq 1$, are isomorphic, and this follows
immediately from the definition of both graphs.
\end{proof}

\begin{figure}[h]
\centering
\input{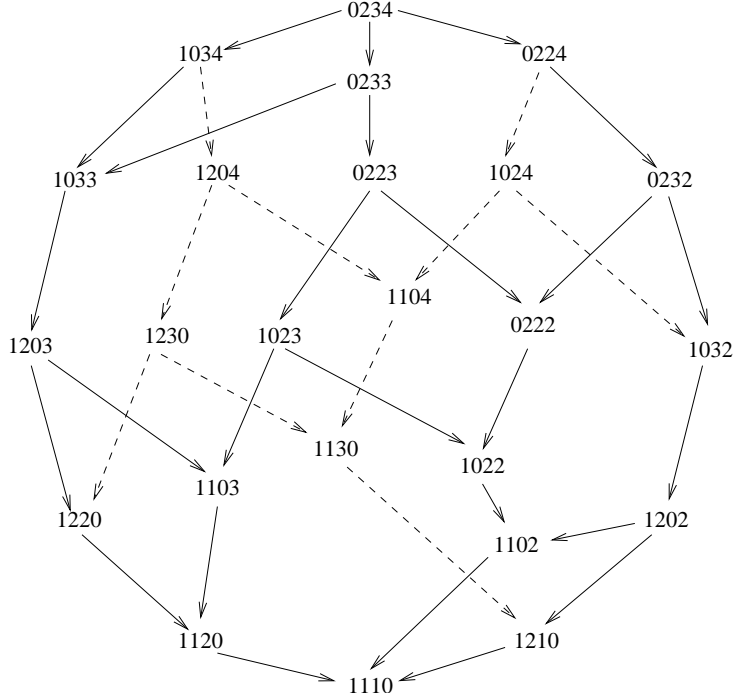}
\caption{The connected component $C(0,2,3,4)$ of $G_{1,3}$}
\label{simplihedre}
\end{figure}

As example, it is shown in Fig.~\ref{simplihedre} the connected
component $C(0,2,3,4)$ of the graph $G_{1,3}$, whose in-vertex is
$(1,1,1,0)$. Notice that it coincides with the 1-skeleton of the
3-dimensional permutohedron $P_3$, which we recall it is obtained
from an octahedron by cutting out 6 small octahedra about its six
vertices. Similarly, the connected components of $G_{1,2}$ (cf.
Fig.~\ref{graf_G12}) were equal to the 1-skeleton of the
2-dimensional permutohedron $P_2$ \footnote{I am very grateful to
the referee for pointing out to me the permutohedron nature of the
connected components of the graphs $G_{s,k}$.}$^)$. Using the
coherence Theorem~\ref{teorema_coherencia}, it is shown below that
this is always true (see Corollary~\ref{corolari}).

Let us now prove the coherence theorem. We proceed in a way very
similar that that followed by MacLane to prove the classical coherence
theorem for monoidal categories (see \cite{sM63}, \cite{sM98}).

\begin{proof}{\it (of Theorem~\ref{teorema_coherencia})}
 Let $v=(i_k,\ldots,i_0)$, $v'=(i'_k,\ldots,i'_0)$ be two arbitrary
vertices in $G_{s,k}$, corresponding to two objects $f,f'$ in $\Cc_{s,k}$.
We have to see that any two different paths between them in $G_{s,k}$ (if
there exists any path at all) correspond to the same morphism in
$\Cc_{s,k}$. We may assume that both vertices belong to the
same connected component, because otherwise there is nothing to be
shown. Let us denote by $C_{s,k}$ this component and let
$v_{in}=(i^{in}_k,\ldots,i^{in}_0)$ be the corresponding
in-vertex. We clearly have a directed path from each vertex $v,v'$ to $v_{in}$
that we may choose in a canonical way, say by always applying in each step the
expanded coherer $1_{g'}\circ\tau_{ij}^n\circ 1_g$ with the least possible
value of $n$ ($n$ will be called the {\it laterality} of the expanded
coherer). This, together with the fact that $\Cc_{s,k}$ is a groupoid, reduces
the proof of the theorem to see that any two directed paths from an arbitrary
vertex $v$ in
$C_{s,k}$ to the vertex $v_{in}$ define the same
isomorphism in $\Cc_{s,k}$. The proof is by induction on the height of
$v\neq v_{in}$. Let $h(v_{in})=h_0$. Hence,
$h(v)\in\{h_0+1,\ldots,h_0+\frac{1}{2}k(k+1)\}$. If $h(v)=h_0+1$, there is
only one path in $C_{s,k}$ from $v$ to $v_{in}$ (a path of length one) and
there is nothing to be shown. Suppose
$h(v)>h_0+1$. We have to distinguish two cases, according to the rank of
$v$. If $r(v)=1$,
there is again a unique directed edge starting at $v$. After crossing that
edge, we get a new vertex $v'$ whose height is $h(v')=h(v)-1$ and the result
follows by the induction hypothesis. Suppose now that $r(v)>1$. In this
case, there are various directed edges starting at $v$, distinguished by the
laterality of the corresponding expanded coherer. By the induction hypothesis,
any two paths starting with the same directed edge in $v$ will define the same
morphism in $\Cc_{s,k}$, because this common first edge will decrease the
height by a
unit. Thus, it only remains to consider the case of two paths $\gamma$,
$\gamma'$ from $v$ starting with
different edges, of lateralities $n$ and $n'$, with $n\neq n'$. The situation
is depicted in Fig~\ref{prova-coherencia}.

\begin{figure}[h]
\centering
\input{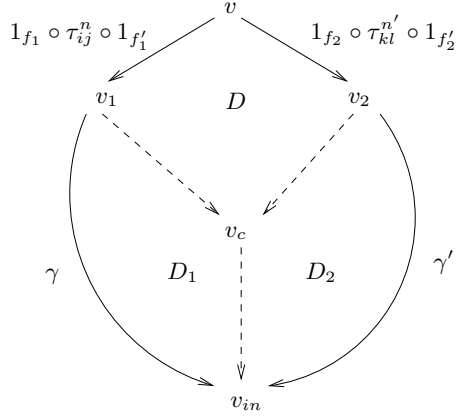}
\caption{Proof of Theorem 4.2}
\label{prova-coherencia}
\end{figure}

It is clear from this figure that we just need to see that both
initial edges can be made to converge to a common vertex $v_c$ in
such a way that the resulting diagram $D$ commutes in $\Cc_{s,k}$,
the corresponding bottom diagrams $D_1$, $D_2$ being commutative by
the induction hypothesis. There are two possibilities, according to
the value of $|n-n'|$. If $|n-n'|=1$, the convergence may be
achieved through an hexagonal diagram, which commutes in
$\Cc_{s,k}$ by the coherence relations. If $|n-n'|>1$, we need just
to apply the expanded coherers with the lateralities interchanged
to get a square which will be commutative in $\Cc_{s,k}$ by the
interchange law.
\end{proof}

Notice that, by the last paragraph in this proof, what it has actually been
shown is that all closed paths in $G_{s,k}$ are the boundary of a union of a
certain number of instances of the hexagonal diagrams giving the
coherence relations (hence, commutative) together with some
quadrilaterals (commutative by the interchange law). Using this,
the above mentioned relation between our graphs $G_{s,k}$ and the
permutohedra easily follows. Let us first recall a few
facts about the permutohedra, defined for the first time by Milgram
\cite{rjM66} (see also
\cite{CM95}, where they are called {\it zilchgons}). For any $k\geq
1$, the {\it permutohedron} $P_k$ is defined as the 
convex hull of the set of points
$(\sigma(1),\ldots,\sigma(k+1))\in\Real^{k+1}$ for all permutations
$\sigma\in S_{k+1}$. It is shown \cite{CM95} that $P_k$ is a
$k$-dimensional convex polyhedron whose $(k+1-r)$-dimensional
faces, for all $r=1,\ldots,k+1$, are indexed by pairs $(p,s)$,
where $p$ is an ordered partition of $\{1,\ldots,k+1\}$, i.e., a
partition of the form
$$
p=\{\{1,\ldots,i_1\},\{i_1+1,\ldots,i_1+i_2\},\ldots
\{i_1+\cdots+i_{r-1}+1,\ldots,i_1+\cdots+i_r\}\}
$$
with $i_1+\cdots+i_r=k+1$ and all $i_j\geq 1$, and $s$ is a shuffle
of type $(i_1,\ldots,i_r)$, namely, a permutation $\sigma\in
S_{k+1}$ such that $\sigma(i)<\sigma(j)$ whenever $i$ and $j$
belong to the same block in the partition. This is equivalent to
label the $(k+1-r)$-dimensional faces by ordered tuples
$(A_1,\ldots,A_r)$ of non-empty disjoint subsets of
$\{1,2,\ldots,k+1\}$ such that $\bigcup_{i=1}^r
A_i=\{1,\ldots,k+1\}$ (the tuple $(A_1,\ldots,A_r)$ corresponding
to a pair $(p,s)$ is obtained by applying the shuffle $s$ to $p$).
In particular, it turns out that the 1-dimensional faces (case
$r=k$) are labelled by pairs $(\sigma,\tau)$, where $\sigma$ is any
permutation in $S_{k+1}$, and $\tau$ is any transposition of the
form $\tau=(i,i+1)$, for some $i\in\{1,\ldots,k\}$; $\tau$ gives,
for the ordering defined by $\sigma$, the two point set in the
corresponding tuple $(A_1,\ldots,A_k)$. For ex., if $k=3$, the
pairs $((123),(34))$, $((14),(12))$ respectively correspond to the
tuples $(\{2\},\{3\},\{1,4\})$ and $(\{4,2\},\{3\},\{1\})$. Such a
pair $(\sigma,(i,i+1))$ represents an edge in $P_k$ between the
vertices
$(\sigma(1),\ldots,\sigma(i),\sigma(i+1),\ldots,\sigma(k+1))$ and
$(\sigma(1),\ldots,\sigma(i+1),\sigma(i),\ldots,\sigma(k+1))$. It
follows that the 1-skeleton of $P_k$ is nothing but the Cayley
graph $Cay(S_{k+1})$ of $S_{k+1}$ with respect to the generators
$\{(12),(23),\ldots,(k,k+1)\}$ (for the definition of the Cayley
graph of a group, see for ex. \cite{CM57}). We then have the
following result, which suggests the name {\it cosemisimplihedra}
for the permutohedra:

\begin{cor} \label{corolari}
For any $k\geq 1$, the connected components of $G_{s,k}$ are
isomorphic to the 1-skeleton of $P_k$.
\end{cor}

\begin{proof}
It is enough to see that the connected component
$C(1,2,\ldots,k+1)$ of $G_{1,k}$, for ex., is isomorphic to
$Cay(S_{k+1})$. If $C(1,\ldots,k+1)_0$ and $Cay(S_{k+1})_0=S_{k+1}$
denote the respective sets of vertices in both graphs, let us
define a map $\Phi:C(1,\ldots,k+1)_0\To Cay(S_{k+1})_0$ by
$$
\Phi(i_0,\ldots,i_k)=(n_1,n_1+1)(n_2,n_2+1)\cdots(n_r,n_r+1),
$$
where $n_1,\ldots,n_r$ are the lateralities of the successive
expanded coherers needed to go from the out-vertex
$(1,2,\ldots,k+1)$ to $(i_0,\ldots,i_k)$. Although there are
can be several paths in $G_{s,k}$ from one vertex to the other, the
corresponding permutation is uniquely defined. Indeed, according to
the remark after the proof of the coherence theorem, any two such
paths are joined through some hexagonal and/or quadrilateral faces.
Now, the hexagonal faces just correspond to the relation
$$
(i,i+1)(i+1,i+2)(i,i+1)=(i+1,i+2)(i,i+1)(i+1,i+2),\quad
i=1,\ldots,k-1
$$
in the symmetric group, while the quadrilaterals correspond to the
relation
$$
(i,i+1)(j,j+1)=(j,j+1)(i,i+1),\quad |i-j|\geq 2
$$
Furthermore, this map is injective, because if two vertices in
$C(1,2,\ldots,k+1)_0$ are mapped to the same permutation, the two
formally different decompositions of the permutation must be
related through the previous relations. But this means that the
corresponding paths in $C(1,\ldots,k+1)$ must be related by
hexagonal and quadrilateral faces as before, so that they
necessarily define a closed path, both final vertices being
equal. Since both sets of vertices have the same cardinal, it 
follows that it is a bijection, and it clearly preserves the edges.
\end{proof}

To finish this section, it is worth emphasizing that, contrary to
what it might seem at first sight, our definition of
2-cosemisimplicial object in $\Cgg$ is not completely equivalent to
a pseudofunctor $\F:\Delta_s\To\Cgg$, where $\Delta_s$ is the
semisimplicial category viewed as a 2-category with only the
identity 2-morphisms
\footnote{This observation has been motivated by a comment of the referee.}$^)$. It is
known that a pseudofunctor $\F:\Cgg\To\Dgg$, with $\Cgg$ and $\Dgg$
2-categories, is equivalent to a 2-functor $F:H(\Cgg)\To\Dgg$,
where $H(\Cgg)$ is a suitable 2-category which depends on $\Cgg$
but not on $\F$. More precisely, it turns out that the inclusion
functor {\bf 2-Cat}$\hookrightarrow ${\bf 2-Cat}$_{ps}$, where {\bf
2-Cat}$_{ps}$ and {\bf 2-Cat} are the categories with objects all
(small) 2-categories and morphisms all pseudofunctors or all
2-functors, respectively, has a left adjoint H:{\bf
2-Cat}$_{ps}\To${\bf 2-Cat} (cf. \cite{sL02}, Prop.4.2
\footnote{I acknowledge the referee for calling my attention to this result,
which seems to be well known among 2-category specialists.}$^)$).
The 2-category $H(\Cgg)$ is in some sense obtained from $\Cgg$ by
making it free with respect to 1-morphisms. Explicitly, it has
as objects the same as $\Cgg$, as 1-morphisms all finite sequences
of composable 1-morphisms in $\Cgg$ (included the empty sequence if
both the domain and codomain objects coincide) and as 2-morphisms
between two such paths all 2-morphisms in $\Cgg$ between the
composite 1-morphisms defined by each path (the composite
1-morphism being the identity when the path is the empty sequence).
Composition of 1-morphisms is given by concatenation and
compositions of 2-morphisms by those in $\Cgg$ in the obvious way.
If $\Cgg=\Delta_s$, the corresponding $H(\Delta_s)$, which will be
denoted by $2\Delta_s$, is a 2-category where, for any two
1-morphisms, we still have at most one 2-morphism between them
(actually, a 2-isomorphism), but they are no longer identities all
of them. For example, for any pair $i,j$ such that $0\leq i<j\leq
n+1$ and any $n\geq 1$, there is in $2\Delta_s$ a (nonidentity)
2-isomorphism
$\beta_{ij}^n:(\epsilon_n^i,\epsilon_{n+1}^j)\Rightarrow(\epsilon_n^{j-1},\epsilon_{n+1}^i)$,
where the $\epsilon_n^i:[n-1]\hookrightarrow[n]$ denote the usual
face morphisms in $\Delta_s$. In $2\Delta_s$, we have a
sub-2-category $2\Delta_s^0$ with the same objects as $2\Delta_s$,
with $2\Delta_s^0([n-1],[n])=2\Delta_s([n-1],[n])$ for all $n\geq
1$, but with $2\Delta_s^0([n-1],[n+k])$, for all $n,k\geq 1$, equal
to the full subcategory of $2\Delta_s([n-1],[n+k])$ whose objects
are only the sequences of composable face morphisms, i.e.,
sequences of length $k+1$ of the form
$(\epsilon_n^{i_0},\ldots,\epsilon_{n+k}^{i_k})$. Notice that this
sub-2-category is biequivalent to $2\Delta_s$, but not equal. Thus,
in $2\Delta^0_s$ the only 1-morphisms defined by sequences of
length 1 are those given by a face morphism in $\Delta_s$, while in
$2\Delta_s$ we further have all sequences of the form $(f)$, for
$f$ any composite of more than one face morphism. We have then the
following:

\begin{prop}
For any 2-category $\Cgg$, a 2-cosemisimplicial object in $\Cgg$ as
defined in Definition~\ref{objecte_2-cosemisimplicial} is
equivalent to a 2-functor $F:2\Delta_s^0\To\Cgg$.
\end{prop}
\begin{proof}
For any 2-category $\Dgg$, a 2-functor $F:\Dgg\To\Cgg$ is
completely defined by the images of the 2-morphisms in any pair of
families of the following type: (1) a family $\Aa_{\Dgg}$ of
``generating'' 2-morphisms in $\Dgg$, by which we mean nonidentity
2-morphisms $\{\tau_{\lambda}\}_{\lambda}$ in $\Dgg$ such that any
nonidentity 2-morphism in $\Dgg$ can be obtained as a (non
necessarily unique) pasting of the $\tau_{\lambda}$ and identity
2-morphisms, and (2) a family $\Bb_{\Dgg}$ of identity 2-morphisms
$\{ 1_{f_{\mu}}\}_{\mu}$ such that any 1-morphism in $\Dgg$ can be
obtained as a (non necessarily unique) composition of the
$f_{\mu}$. This is because a 2-functor preserves vertical and
horizontal compositions (hence, also pastings) and the identity
2-morphisms, together with the fact that $1_{g\circ f}=1_g\circ
1_f$. Furthermore, the images of these 2-morphisms can be chosen
arbitrarily except that all possible relations
between them have to be preserved. If $\Dgg=2\Delta_s^0$, a pair of
families as above is given by
\begin{align*}
\Aa_{2\Delta_s^0}&=\{\beta_{ij}^n,\ 0\leq i<j\leq n+1,\ n\geq 1\}
\\ \Bb_{2\Delta_s^0}&=\{1_{\epsilon_n^i},\ i=0,1,\ldots,n,\ n\geq 1\}
\end{align*}
Indeed, the face morphisms generate all morphisms in $\Delta_s$
and, given two sequences
$(\epsilon_n^{i_0},\ldots,\epsilon_{n+k}^{i_k})$,
$(\epsilon_n^{i'_0},\ldots,\epsilon_{n+k}^{i'_k})$ defining the
same composite morphism in $\Delta_s$, they can be connected by a
pasting of the $\beta_{ij}^r$ to the common canonical decomposition
$(\epsilon_n^{j_0},\ldots,\epsilon_{n+k}^{j_k})$ with
$j_0<j_1<\cdots<j_k$. Hence, the corresponding (unique) 2-morphism
between both sequences is really a pasting of the $\beta_{ij}^r$.
It also follows from the uniqueness of the 2-morphisms in
$2\Delta_s^0$ that the $\beta_{ij}^n$ satisfy the relations
\begin{align*}
(1_{\epsilon_{n+2}^i}\circ\beta_{j-1,k-1}^n)\cdot(\beta_{ik}^{n+1}&\circ
1_{\epsilon_n^{j-1}})\cdot(1_{\epsilon_{n+2}^k}\circ\beta_{ij}^n)=
\\ &=(\beta_{ij}^{n+1}\circ
1_{\epsilon_n^{k-2}})\cdot(1_{\epsilon_{n+2}^j}\circ\beta_{i,k-1}^n)\cdot(\beta_{jk}^{n+1}\circ
1_{\epsilon_n^i})
\end{align*}
for all $0\leq i<j<k\leq n+2$ and all $n\geq 1$. Moreover, our
coherence theorem (Theorem~\ref{teorema_coherencia}) implies that
any other relation between the $\beta_{ij}^n$ is a consequence of
these relations. Hence, giving a 2-functor $F:2\Delta_s^0\To\Cgg$
is indeed equivalent to give arbitrary 1-morphisms $\partial_n^i$
in $\Cgg$ (the images of the 2-morphisms $1_{\epsilon_n^i}$) and
2-morphisms $\tau_{ij}^n$ (the images of the 2-morphisms
$\beta_{ij}^n$) in $\Cgg$ satisfying the coherence relations in
Definition~\ref{objecte_2-cosemisimplicial}.
\end{proof}
Such a 2-functor $F:2\Delta_s^0\To\Cgg$, however, will not extend
uniquely to a 2-functor $\widetilde{F}:2\Delta_s\To\Cgg$. Thus,
$2\Delta_s$ contains 2-morphisms 
$(f)\Rightarrow(\epsilon_n^{i_0},\ldots,\epsilon_{n+k}^{i_k})$ with
$f=\epsilon_{n+k}^{i_k}\circ\cdots\circ\epsilon_n^{i_0}$, $k\geq
1$ which are not given by a pasting of the
$\beta_{ij}^n$ and, hence, such that their images are not determined by the
images of the $\beta_{ij}^n$. The reader may easily check that
a right family $\Aa_{2\Delta_s}$ of generating 2-morphisms for
$2\Delta_s$ is given by the family $\Aa_{2\Delta_s^0}$ together with
the 2-morphisms 
$$
\alpha_{i_0,\ldots,i_k}^n:(\epsilon_{n+k}^{i_k}\circ\cdots\circ\epsilon_n^{i_0})\Rightarrow
(\epsilon_n^{i_0},\ldots,\epsilon_{n+k}^{i_k})
$$
for all $0\leq i_0<i_1<\cdots<i_k\leq n+k$ and all $n\geq 1$. To
define a 2-functor $\widetilde{F}:2\Delta_s\To\Cgg$, and hence a
pseudofunctor $\F:\Delta_s\To\Cgg$, we will need to give images
of the 2-morphisms in this additional family satisfying the
appropriate relations. It seems possible, however, that all extensions
$\widetilde{F}:2\Delta_s\To\Cgg$ of 
$F$ turn out to be equivalent in a suitable sense (in the same way
as the extension of a functor defined on the skeleton of a category
to the whole category is unique up to isomorphism), but we did not
explore that any further.


\section{Cochain complexes from 2-cosemisimplicial objects in ${\bf
    Cat}_K$}

Given a cosemisimplicial object in an abelian category, it is usual
to consider the corresponding cochain complex and cohomology.
Hence, the following question naturally raises: what are the
analogs of these cochain complexes and their cohomologies in the
case of a 2-cosemisimplicial object in a 2-category? As in the
categorical setting, it is expected that finding these analogs will
require restricting to suitable {\it abelian 2-categories}, for
which hypothetical {\it 2-cochain complexes} will make sense.
However, we will not pursue  this direction here. Instead, the
purpose of this section is to show that usual cochain complexes of
$K$-modules may still be constructed from certain {\it enhanced}
2-cosemisimplicial objects in a particular 2-category. Namely, the
2-category ${\bf Cat}_K$ having as objects the (small) $K$-linear
categories, as 1-morphisms the $K$-linear functors and as
2-morphisms the natural transformations.
As an example, which was our original motivation, we show in the
next section that the purely pseudofunctorial deformation complex
introduced in \cite{jE1} for any $K$-linear pseudofunctor $\F$ may
be obtained in this way from a suitable enhanced 2-cosemisimplicial
object in ${\bf Cat}_K$ associated to $\F$.

Suppose we are given a 2-cosemisimplicial object $\Cc^{\bullet}$ in
${\bf Cat}_K$ and let $F_n^i:\Cc^{n-1}\rightarrow\Cc^n$
($i=0,1,\ldots,n$, $n\geq 1$) and $\tau_{ij}^n:F_{n+1}^j\circ
F_n^i\Rightarrow F^i_{n+1}\circ F^{j-1}_n$ ($0\leq i<j\leq n+1$,
$n\geq 1$) be the corresponding coface functors and
cosemisimplicial coherers, which are natural isomorphisms in this
case. To simplify notation, we shall write $F^{i_0,\ldots,i_k}_n$
to denote the composite functor $F^{i_k}_{n+k}\circ
F^{i_{k-1}}_{n+k-1}\circ\cdots\circ F^{i_0}_n$ ($n,k\geq 1$).
According to Theorem~\ref{teorema_coherencia}, for all $n,k\geq 1$
and $(i_0,\ldots,i_k)\neq (j_0,\ldots,j_k)$, with
$i_q,j_q\in\{0,1,\ldots,n+q\}$ and $q=0,1,\ldots,k$, there exists
at most one canonical natural isomorphism from
$F^{i_0,\ldots,i_k}_n$ to $F^{j_0,\ldots,j_k}_n$, given by pasting
the appropriate coherers $\tau_{ij}^n$'s and/or its inverses. It
will be denoted by $\tau_{(i_0,\ldots,i_k),(j_0,\ldots,j_k)}^n$.
Notice that such canonical isomorphisms may not exist, depending on
the $(k+1)$-tuples $(i_0,\ldots,i_k)$ and $(j_0,\ldots,j_k)$. This
is because, as seen before, the graph $G_{n,k}$ is not connected.
For example, there is no canonical path between $F_n^{1,1}$ and
$F_n^{0,0}$ nor between $F_n^{1,1,0,3}$ and $F_n^{1,2,3,4}$. When
$(i_0,\ldots,i_k)=(j_0,\ldots,j_k)$, we will agree that
$\tau_{(i_0,\ldots,i_k),(i_0,\ldots,i_k)}^n$ denotes the
corresponding identity natural transformation.

Roughly, the method of getting cochain complexes of $K$-modules
from the 2-cosemisimplicial object $\Cc^{\bullet}$ consists of the
following. For all $n\geq 0$, choose a pair of objects $X_n$,
$X'_n$ in each category $\Cc^n$, take for each such pair the
corresponding $K$-modules of morphisms ${\rm
Hom}_{\Cc^n}(X_n,X'_n)$ (they are indeed $K$-modules because
$\Cc^n$ is $K$-linear) and define coboundary maps between them
using the coface functors $F_n^i$, which are $K$-linear. More
explicitly, we would like these coboundary maps $\delta:{\rm
  Hom}_{\Cc^{n-1}}(X_{n-1},X'_{n-1})\rightarrow {\rm
  Hom}_{\Cc^n}(X_n,X'_n)$ to be of the form
\begin{equation} \label{intent_delta}
\delta(\varphi)\approx \sum_{i=0}^n (-1)^iF_n^i(\varphi)
\end{equation}
for all $\varphi\in{\rm Hom}_{\Cc^{n-1}}(X_{n-1},X'_{n-1})$. This
procedure, however, makes no sense in general, because the
$F_n^i(\varphi)$ belong to different $K$-modules of
morphisms for different
values of $i\in\{0,1,\ldots,n\}$ (they have different domains and
codomains). This could be
easily overcomed if all such domains and codomains were (canonically)
isomorphic to the
corresponding reference objects $X_n$ and $X'_n$,
respectively, because we can then get morphisms in ${\rm
  Hom}_{\Cc^n}(X_n,X'_n)$ by just taking the composite of each term
$F_n^i(\varphi)$ with the appropriate (canonical) isomorphisms on the left and
on the right. However, this will not be true for randomly chosen objects $X_n$ and $X_n'$. One may try to fix that by choosing an object
$X\in|\Cc^0|$ and taking $X_n$ and $X'_n$, for all $n\geq 1$, equal to some
iterated images of $X$ by the coface functors $F_n^i$. For
example, for $n\geq 1$, we could inductively define
\begin{align}
X_n&:=F_n^n(X_{n-1}) \label{objectes_referencia_1} \\
X'_n&:=F_n^{n-1}(X'_{n-1}) \label{objectes_referencia_2}
\end{align}
with $X_0=X'_0=X$. In this way, both the domain and codomain of
$F_n^i(\varphi)$, for all $i=0,\ldots,n$, will be of the form
$F_1^{i_0,\ldots,i_{n-1}}$ for some
$n$-tuples of positive integers $(i_0,\ldots,i_{n-1})$, so that they can be
related via the natural isomorphisms $\tau_{ij}^n$. Even in this way,
however, the problem turns out to persist because of the
non-connectedness of the graphs $G_{1,n-1}$. Actually, the problem persists
independently of how
the references $X_n$ and $X'_n$ are chosen among all
possible iterated images
of $X$. This is easily seen by considering the
cases $n=1$ and $n=2$. Suppose we take $X_1=F_1^1(X)$. Then, for any
$\varphi:X_1\rightarrow X'_1$, the domains of $F_2^0(\varphi)$,
$F_2^1(\varphi)$ and $F_2^2(\varphi)$ will respectively be $F_1^{1,0}(X)$,
$F_1^{1,1}(X)$ and  $F_1^{1,2}(X)$. But a glance to the graph $G_{1,1}$
immediately shows that there is no choice for $X_2=F_2^i(X_1)$ such that it is
simultaneously canonically isomorphic to these three domains.

The above discussion shows that to define cochain complexes by this
method, with the coboundary maps given by
Equation~\ref{intent_delta}, we need some additional hypothesis on
the 2-cosemisimplicial object $\Cc^{\bullet}$. This leads us to
introduce the following definition.

\begin{defn}
Let $\Cgg$ be any 2-category. By an {\sl enhanced}
2-cosemisimplicial object in $\Cgg$ we shall mean a
2-cosemisimplicial object $(X^{\bullet},\partial,\tau)$ in $\Cgg$
together with a 2-isomorphism
$\phi:\partial_1^1\Rightarrow\partial_1^0$ such that
\begin{equation} \label{relacio_coherencia_phi}
\tau^1_{0,1}\cdot(1_{\partial_2^1}\circ\phi)\cdot\tau^1_{1,2}=(1_{\partial_2^0}\circ\phi)\cdot\tau^1_{0,2}\cdot(1_{\partial_2^2}\circ\phi)
\end{equation}
\end{defn}
As the coherence relations on the $\tau_{ij}^n$'s, the above condition on
 $\phi$ is related to a coherence theorem. To state this theorem, let us
 denote by $\Cc^{\phi}_{1,k}$, for all $k\geq 1$, the subcategory of
 $\Cgg(X^0,X^{k+1})$ with objects
 the same as in $\Cc_{1,k}$ (namely, the $\partial$-paths), but whose
 morphisms are all possible composites of expanded coherers of
 $X^{\bullet}$ and expansions of $\phi$ (i.e., 2-isomorphisms of the form
 $1_f\circ\phi:f\circ\partial_1^1\Rightarrow f\circ\partial_1^0$ for some
 $\partial$-path $f$). The new  coherence theorem states then the following:

\begin{thm} \label{segon_teorema_coherencia}
Let $k\geq 1$. Then, for any two objects $f,f'$ in $\Cc^{\phi}_{1,k}$, there
is one and only one morphism (actually an isomorphism) in $\Cc^{\phi}_{1,k}$
from $f$ to $f'$.
\end{thm}
\begin{proof}
Let $G_{1,k}^{\phi}$ be the graph with vertices all $\partial$-paths
 $f:X^0\rightarrow X^{k+1}$ and edges all the expanded coherers and
 expansions of $\phi$. In particular, $G_{1,k}^{\phi}$ contains $G_{1,k}$ as a
 subgraph (see Figure~\ref{graf_G12_phi} for the case $k=2$).
\begin{figure}[h]
\centering
\input{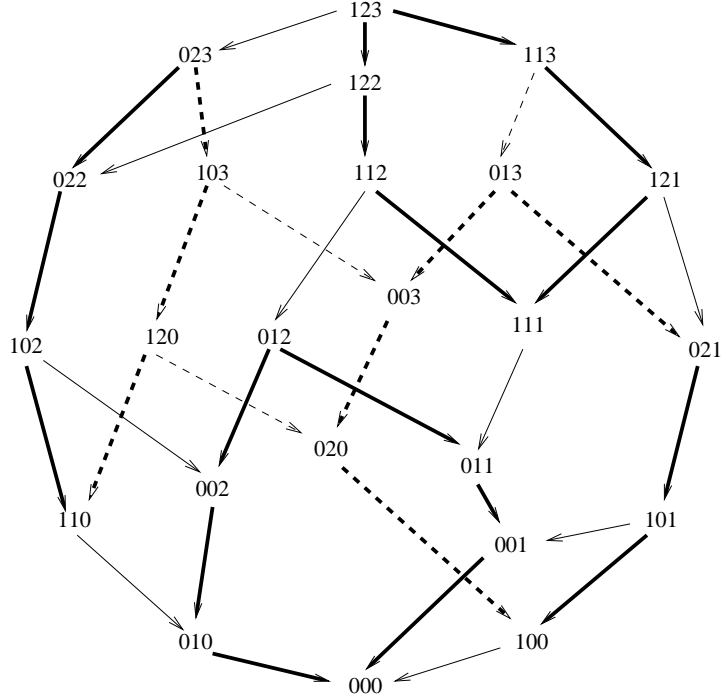}
\caption{The graph $G^{\phi}_{1,2}$ (the edges in the four hexagons defining $G_{1,2}$ are drawn in bold
  solid or dashed arrows to distinguish them from the additional edges
  corresponding to the expansions of $\phi$).}
\label{graf_G12_phi}
\end{figure}
As in the previous
 section, we identify a $\partial$-path $f$ with the corresponding
 $(k+1)$-tuple $(i_0,\ldots,i_k)$. Let us first prove that, given two
 arbitrary vertices $(i_0,\ldots,i_k)$ and $(i'_0,\ldots,i'_k)$, there always
 exist a path in $G_{1,k}^{\phi}$ between them. Clearly, it is enough to
 prove the assertion in the special case when both vertices are
 in-vertices. Otherwise, one immediately obtains the desired path by
 connecting each vertex to the corresponding in-vertex and
 adding any path between both in-vertices. To prove
 the claim for two in-vertices, observe that all in-vertices in
 $G_{1,k}^{\phi}$ are of the form $(1,\ldots,
 1,0,\ldots,0)$, the number of 1's plus the number of 0's being equal to
 $k+1$. Starting at any such in-vertex and via the appropriate expansion of
 $\phi$, we can move to the neighbour ``dual''
 vertex $(0,1,\ldots,1,0,\ldots,0)$, differing from it just in the first
 component. This is not an in-vertex, but it can be connected to the
 corresponding in-vertex through a path of expanded coherers. This new
 in-vertex will have one more zero than the initial one. Iterating
 this process, one finally gets the in-vertex
 $(0,\ldots,0)$. Since this may be done for any initial in-vertex, we conclude
 that two arbitrary in-vertices are indeed connected in $G_{1,k}^{\phi}$. Let
 us now prove that all paths in $G_{1,k}^{\phi}$ between two arbitrary
 vertices $(i_0,\ldots,i_k)$ and $(i'_0,\ldots,i'_k)$ define the same morphism
 in $\Cc_{1,k}^{\phi}$. Since $\Cc_{1,k}^{\phi}$ is also a groupoid, it
 suffices to prove the assertion when
 $(i'_0,\ldots,i'_k)=(0,\ldots,0)$. We proceed again by induction on the
 height of $v=(i_0,\ldots,i_k)$. If $h(v)=0$, $v$ necessarily coincides with
 $(0,\ldots,0)$ and there is nothing to be shown. Suppose then that $h(v)\geq
 1$ and let $\gamma$, $\gamma'$ be two directed paths starting at $v$. If the
 first edges in both $\gamma$ and $\gamma'$ coincide, the
 result follows directly by induction. Otherwise, the argument is similar to
 that used in the proof of Theorem~ \ref{teorema_coherencia}. Namely, we see
 that both initial edges can be made to converge to a common vertex $v_c$ in
 such a way that the resulting diagram $D$ commutes in $\Cc_{1,k}^{\phi}$. We
 have to distinguish three possibilities:

\noindent{(i)} Both first edges
 are two different expanded coherers: in this case, the convergence is achieved
 in exactly the same way as in the proof of Theorem~\ref{teorema_coherencia}.

\noindent{(i)} One edge is an expansion $1_f\circ\phi$ of $\phi$ while the other one is
 an expanded coherer $1_{f'}\circ\tau_{ij}^s\circ 1_{f''}$ with laterality
 $s\geq 2$: in this case the commutative diagram $D$ is the square defined by
 the equality
$$
(1_{f'}\circ\tau_{ij}^s\circ 1_{\overline{f}''})\cdot(1_f\circ\phi)=
(1_{\overline{f}}\circ\phi)\cdot(1_{f'}\circ\tau_{ij}^s\circ 1_{f''})
$$
which holds by the interchange law.

\noindent{(i)}One edge is an expansion $1_f\circ\phi$ of $\phi$ while the other one is
 an expanded coherer of the form $1_{f'}\circ\tau_{ij}^1$, with laterality
 $s=1$: in this case the commutative
 diagram $D$ is just that defined by Equation~\ref{relacio_coherencia_phi},
 which holds by hypothesis.
\end{proof}
These unique isomorphisms between the objects in
$\Cc_{1,k}^{\phi}$ will be called the {\it canonical enhanced
  2-isomorphisms} and denoted by
$\tau^{\phi}_{(i_0,\ldots,i_k),(j_0,\ldots,j_k)}$. Notice that, when the pair
$(i_0,\ldots,i_k),(j_0,\ldots,j_k)$ is such that there already exists a
path of expanded coherers in $G_{1,k}$ between the corresponding
$\partial$-paths, this canonical enhanced 2-isomorphism
$\tau^{\phi}_{(i_0,\ldots,i_k),(j_0,\ldots,j_k)}$ coincides with the
canonical 2-isomorphism $\tau^1_{(i_0,\ldots,i_k),(j_0,\ldots,j_k)}$ defined
in the previous section.

\begin{rem}
{\rm Suprisingly, the graph $G_{1,2}^{\phi}$ turns out to be
  isomorphic to the connected components of $G_{1,3}$ and, hence, to the
  1-skeleton of the 3-dimensional permutohedron $P_3$ (cf. Figs.~\ref{simplihedre}
  and \ref{graf_G12_phi}). This suggests that the same may be
  true for all $k\geq 2$.}
\end{rem}

We may now carry out the above program. Let
$(\Cc^{\bullet},F,\tau,\phi)$ be an enhanced 2-cosemisimplicial
object in ${\bf Cat}_K$ and let us fix an object $X\in|\Cc^0|$. For
all $n\geq 1$, choose once and for all $n$-tuples of nonnegative
integers $(\mu^n_1,\ldots,\mu^n_n)$ and $(\nu^n_1,\ldots,\nu^n_n)$,
with $\mu^n_q,\nu^n_q\in\{0,\ldots,q\}$, and define objects
$X_n,X'_n\in|\Cc^n|$ by
\begin{align*}
X_n&=F_1^{\mu^n_1,\ldots,\mu^n_n}(X) \\ X'_n&=F_1^{\nu^n_1,\ldots,\nu^n_n}(X)
\end{align*}
They will
be called the {\it domain} and {\it codomain
  reference objects} in $\Cc^n$, respectively. According to the previous
theorem, for all $n\geq 1$ and $i=0,1,\ldots,n$, we have canonical enhanced
2-isomorphisms
\begin{align*}
\tau^{\phi}_{(\nu^{n-1}_1,\ldots,\nu^{n-1}_{n-1},i),(\nu^n_1,\ldots,\nu^n_{n-1},\nu^n_n)}&:F_1^{\nu^{n-1}_1,\ldots,\nu^{n-1}_{n-1},i}\Rightarrow
F_1^{\nu^n_1,\ldots,\nu^n_{n-1},\nu^n_n} \\
\tau^{\phi}_{(\mu^n_1,\ldots,\mu^n_{n-1},\mu^n_n),(\mu^{n-1}_1,\ldots,\mu^{n-1}_{n-1},i)}&:F_1^{\mu^n_1,\ldots,\mu^n_{n-1},\mu^n_n}\Rightarrow
F_1^{\mu^{n-1}_1,\ldots,\mu^{n-1}_{n-1},i}
\end{align*}
Hence, by taking the corresponding $X$-components, we get isomorphisms
\begin{align*}
\alpha^{\phi}_{i,n}\equiv\left(\tau^{\phi}_{(\nu^{n-1}_1,\ldots,\nu^{n-1}_{n-1},i),(\nu^n_1,\ldots,\nu^n_{n-1},\nu^n_n)}\right)_X&:F_n^i(X'_{n-1})\rightarrow
X_n' \\
\beta^{\phi}_{i,n}\equiv\left(\tau^{\phi}_{(\mu^n_1,\ldots,\mu^n_{n-1},\mu^n_n),(\mu^{n-1}_1,\ldots,\mu^{n-1}_{n-1},i)}\right)_X&:X_n\rightarrow F_n^i(X_{n-1})
\end{align*}
Let us further denote by
$M^n$, for all $n\geq 1$, the $K$-module of morphisms
$$
M^n:={\rm Hom}_{\Cc^n}(X_n,X'_n)
$$
We have then the following:

\begin{thm} \label{complex_derivat}
The above $K$-modules
$M^1,M^2,\ldots$ together with the coboundary maps $\delta:M^{n-1}\rightarrow
M^n$, $n\geq 2$, given by
\begin{equation} \label{operadors_covora}
\delta(\varphi)=\sum_{i=0}^n (-1)^i\ \alpha^{\phi}_{i,n}\circ
F_n^i(\varphi)\circ\beta^{\phi}_{i,n}
\end{equation}
define a cochain complex. Furthermore, the cochain complexes defined in this
way by different choices of the reference objects $X_n$, $X'_n$, $n\geq 1$,
and by different objects $Y\in|\Cc^0|$ isomorphic to $X$ are all isomorphic.
\end{thm}
\begin{proof}
Since the functors $F_n^i$ are $K$-linear, the coboundary maps are indeed
$K$-linear. To see that $\delta^2=0$, notice first that, by the naturality of
the $\tau_{ij}^n$'s, we have
$$
(F_{n+1}^j\circ
F_n^i)(\varphi)=(\tau_{ij}^n)^{-1}_{X'_{n-1}}\circ(F_{n+1}^i\circ
F_n^{j-1})(\varphi)\circ(\tau_{ij}^n)_{X_{n-1}}
$$
for any $\varphi:X_{n-1}\rightarrow X'_{n-1}$ and all $0\leq i<j\leq
n+1$. Then, proceeding in the usual way, we have
\begin{align*}
\delta^2(\varphi)&=\sum_{0\leq i<j\leq n+1}(-1)^{i+j}\alpha^{\phi}_{j,n+1}\circ F_{n+1}^j(\alpha^{\phi}_{i,n})\circ
(\tau^n_{ij})^{-1}_{X'_{n-1}} \\ &\hspace{3 truecm}\circ(F_{n+1}^i\circ F_n^{j-1})(\varphi)\circ(\tau^n_{ij})_{X_{n-1}}\circ
F_{n+1}^j(\beta^{\phi}_{i,n})\circ\beta^{\phi}_{j,n+1} \\
&\ +\sum_{n\geq i\geq j\geq 0}(-1)^{i+j}\alpha^{\phi}_{j,n+1}\circ F_{n+1}^j(\alpha^{\phi}_{i,n})\circ
(F_{n+1}^j\circ F_n^i)(\varphi)\circ
F_{n+1}^j(\beta^{\phi}_{i,n})\circ\beta^{\phi}_{j,n+1} \\
&=\sum_{0\leq j\leq i\leq n}(-1)^{i+j+1}\alpha^{\phi}_{i+1,n+1}\circ F_{n+1}^{i+1}(\alpha^{\phi}_{j,n})\circ
(\tau^n_{j,i+1})^{-1}_{X'_{n-1}} \\ &\hspace{3 truecm}\circ(F_{n+1}^j\circ F_n^i)(\varphi)\circ(\tau^n_{j,i+1})_{X_{n-1}}\circ
F_{n+1}^{i+1}(\beta^{\phi}_{j,n})\circ\beta^{\phi}_{i+1,n+1} \\
&\ +\sum_{n\geq i\geq j\geq 0}(-1)^{i+j}\alpha^{\phi}_{j,n+1}\circ F_{n+1}^j(\alpha^{\phi}_{i,n})\circ
(F_{n+1}^j\circ F_n^i)(\varphi)\circ
F_{n+1}^j(\beta^{\phi}_{i,n})\circ\beta^{\phi}_{j,n+1}
\end{align*}
the last equality being obtained by a suitable reindexation in the first
sum. Hence, the proof reduces to see that the
$\alpha$'s, $\beta$'s and $\tau$'s satisfy the equations
\begin{align}
\alpha^{\phi}_{j,n+1}\circ
F_{n+1}^j(\alpha^{\phi}_{i,n})\circ(\tau^n_{j,i+1})_{X'_{n-1}}&=\alpha^{\phi}_{i+1,n+1}\circ
F_{n+1}^{i+1}(\alpha^{\phi}_{j,n}) \label{condicio1_alfa_beta} \\ (\tau^n_{j,i+1})_{X_{n-1}}\circ F^{i+1}_{n+1}(\beta^{\phi}_{j,n})\circ\beta^{\phi}_{i+1,n+1}&=F^j_{n+1}(\beta^{\phi}_{i,n})\circ\beta^{\phi}_{j,n+1}
\end{align}
for all $0\leq j\leq i\leq n$ ($n\geq 2$). Now, from the very definition of
all the involved terms, we have that the left-hand side in the first equality
is nothing but the $X$-component of the canonical enhanced 2-isomorphism
\begin{align*}
&\tau^{\phi}_{(\nu^n_1,\ldots,\nu^n_n,j),(\nu^{n+1}_1,\ldots,\nu^{n+1}_n,\nu^{n+1}_{n+1})}\cdot\left(1_{F_{n+1}^j}\circ\tau^{\phi}_{(\nu^{n-1}_1,\ldots,\nu^{n-1}_{n-1},i),(\nu^n_1,\ldots,\nu^n_{n-1},\nu^n_n)}\right)
  \\ &\hspace{6 truecm}\cdot\left(\tau_{j,i+1}^n\circ
  1_{F_1^{\nu^{n-1}_1,\ldots,\nu^{n-1}_{n-1}}}\right)
\end{align*}
while the right-hand side is the $X$-component of the canonical enhanced
2-isomorphism
$$
\tau^{\phi}_{(\nu^n_1,\ldots,\nu^n_n,i+1),(\nu^{n+1}_1,\ldots,\nu^{n+1}_n,\nu^{n+1}_{n+1})}\cdot\left(1_{F_{n+1}^{i+1}}\circ\tau^{\phi}_{(\nu^{n-1}_1,\ldots,\nu^{n-1}_{n-1},j),(\nu^n_1,\ldots,\nu^n_{n-1},\nu^n_n)}\right)
$$
By the coherence
Theorem~\ref{segon_teorema_coherencia}, both 2-isomorphisms coincide. The
second equality is shown in a similar way. Let us now prove that the
isomorphism class of the cochain complex is independent of the chosen
references. Indeed, suppose
we choose other references $\overline{X}_n$, $\overline{X}'_n$, defined by
$n$-tuples $(\overline{\mu}^n_1,\ldots,\overline{\mu}^n_n)$ and
$(\overline{\nu}^n_1,\ldots,\overline{\nu}^n_n)$, for all $n\geq 1$. Then, the
new $K$-module
$$
\overline{M}^n={\rm Hom}_{\Cc^n}(\overline{X}_n,\overline{X}'_n)
$$
$n\geq 1$, is isomorphic to the old one $M^n$ through the
isomorphism $f^n:M^n\rightarrow\overline{M}^n$ defined by
$$
f^n(\varphi)=\left(\tau^{\phi}_{(\nu^n_1,\ldots,\nu^n_n),(\overline{\nu}^n_1,\ldots,\overline{\nu}^n_n)}\right)_X\circ\varphi\circ\left(\tau^{\phi}_{(\overline{\mu}^n_1,\ldots,\overline{\mu}^n_n),(\mu^n_1,\ldots,\mu^n_n)}\right)_X
$$
The coboundary operators
$\overline{\delta}:\overline{M}^{n-1}\rightarrow\overline{M}^n$ are defined as
before, except that we have to use now the isomorphisms
$\overline{\alpha}^{\phi}_{i,n}$ and $\overline{\beta}^{\phi}_{i,n}$
corresponding to the new references. It easily follows again from
Theorem~\ref{segon_teorema_coherencia} that the $f^n$ define a cochain
map. Finally, suppose we choose another object $Y\cong X$, $Y\in|\Cc^0|$ and
let us denote by $N^n$ the corresponding $K$-modules, namely, for all $n\geq
1$,
$$
N^n={\rm Hom}_{\Cc^n}(Y_n,Y'_n)
$$
where
\begin{align*}
Y_n&=F_1^{\mu^n_1,\ldots,\mu^n_n}(Y) \\ Y'_n&=F_1^{\nu^n_1,\ldots,\nu^n_n}(Y)
\end{align*}
The coboundary maps are defined as before but using the $Y$-component of the
 corresponding canonical enhanced 2-isomorphisms, i.e., the isomorphisms
\begin{align*}
\gamma^{\phi}_{i,n}\equiv\left(\tau^{\phi}_{(\nu^{n-1}_1,\ldots,\nu^{n-1}_{n-1},i),(\nu^n_1,\ldots,\nu^n_{n-1},\nu^n_n)}\right)_Y&:F_n^i(Y'_{n-1})\rightarrow
Y_n' \\
\eta^{\phi}_{i,n}\equiv\left(\tau^{\phi}_{(\mu^n_1,\ldots,\mu^n_{n-1},\mu^n_n),(\mu^{n-1}_1,\ldots,\mu^{n-1}_{n-1},i)}\right)_Y&:Y_n\rightarrow F_n^i(Y_{n-1})
\end{align*}
instead of the $\alpha^{\phi}_{i,n}$ and $\beta^{\phi}_{i,n}$. Now, if
$h:X\rightarrow Y$ is an isomorphism, it follows immedatiely from the
naturality of the canonical enhanced 2-isomorphisms that
$$
F_1^{\nu_1^n,\ldots,\nu_{n-1}^{n-1}\nu_n^n}(h)\circ\alpha^{\phi}_{i,n}=\gamma^{\phi}_{i,n}\circ
F_1^{\nu_1^{n-1},\ldots,\nu_{n-1}^{n-1},i}(h)
$$
and that a similar relation holds between the $\beta^{\phi}_{i,n}$ and the
$\eta^{\phi}_{i,n}$. Then, defining isomorphisms
$f^n:M^n\rightarrow N^n$ by
$$
f^n(\varphi)=F_1^{\nu_1^n,\ldots,\nu_n^n}(h)\circ\varphi\circ
F_1^{\mu_1^n,\ldots,\mu_n^n}(h^{-1})
$$
it is easily checked that we obtain an isomorphism of cochain complexes.
\end{proof}

\begin{rem}
{\rm Enhanced 2-cosemisimplicial objects are needed to define
cochain complexes with coboundary maps of the form
(\ref{intent_delta}), where the alternating sum is over all coface
functors $F_n^i$, for all $i=0,\ldots,n$. However, it is well-known
that, given a cosemisimplicial object in an abelian category, there
are other cochain complexes that may be defined from it. For
example, one may define the so-called {\sl path space cochain
complex} (see \cite{cW94}), a cochain complex starting at $X^1$
instead of at $X^0$ and whose coboundary maps are given by the
alternating sum
$\delta=\delta^1_n-\delta^2_n+\cdots+(-1)^{n+1}\delta^n_n$, where
the first coface map $\delta^0_n$ has been omitted. In this sense,
it is worth to point out that some of these alternative cochain
complexes can be defined even for {\sl arbitrary}
2-cosemisimplicial objects in ${\bf Cat}_K$. In particular, this is
the case for the dual path space cochain complex of the previous
path space, which is a cochain complex starting at $X^2$ and with
coboundary maps given by
$\delta=\delta^1_n-\delta^3_n+\cdots+(-1)^n\delta^{n-1}_n$ (both
$\delta^0_n$ and $\delta^n_n$ are omitted). We leave to the reader
to check that it is indeed possible to choose reference objects
$X_n$, $X_n'$ in such a way that all the involved domains and
codomains of the maps $F_n^i$, for all $i=1,\ldots,n-1$, belong to
the same connected component of the graph $G_{1,n-1}$, so that no
enhancement is needed in this case to construct a cochain complex
by the previous method.}
\end{rem}


\section{2-cosemisimplicial object of a pseudofunctor and the deformation
  complex}

We are now in a position that enables us to prove the result
mentioned in the introduction. Namely, that associated to any
pseudofunctor $\F$ there is a 2-cosemisimplicial object in {\bf
Cat} and that, when $\F$ is $K$-linear, the cochain complex
$X^{\bullet}(\F)$ introduced in \cite{jE1} is the cochain complex
obtained by the above method from the corresponding
2-cosemisimplicial object in ${\bf Cat}_K$.

Let $\F:\Cgg\rightarrow\Dgg$ be an arbitrary pseudofunctor between
2-categories. Included in these data, we have three collections of
functors. Namely, the composition functors of $\Cgg$ and $\Dgg$
\begin{align*}
c^{\Cgg}_{X,Y,Z}&:\Cgg(X,Y)\times\Cgg(Y,Z)\To\Cgg(X,Z),\quad  X,Y,Z\in|\Cgg|
\\  c^{\Dgg}_{U,V,W}&:\Dgg(U,V)\times\Dgg(V,W)\To\Dgg(U,W),\quad U,V,W\in|\Dgg|
\end{align*}
and the functors
$$
\F_{X,Y}:\Cgg(X,Y)\To\Dgg(\F(X),\F(Y)),\quad  X,Y\in|\Cgg|
$$
defining the action of $\F$ on the 1- and 2-morphisms. From such
functors, and given $X_0,\ldots,X_n\in|\Cgg|$, we may construct various
iterates, differing in the way they apply an arbitrary
path of 1-morphisms in $\Cgg$
$$
\gamma:X_0\stackrel{f_1}{\To} X_1\To\cdots\To X_{n-1}\stackrel{f_n}{\To}
  X_n
$$
to a path in $\Dgg$. More precisely, we define the
following of iterate of $\F$.
\begin{defn}
Given $n\geq 1$ and $X_0,\ldots,X_n\in|\Cgg|$, an
$\F_{X_0,\ldots,X_n}$-{\it iterate} is any functor
$$
H_{X_0,\ldots,X_n}:\Cgg(X_0,X_1)\times\cdots\times\Cgg(X_{n-1},X_n)\To\Dgg(\F(X_0),\F(X_n))
$$
obtained as a composite of products of the functors
$\F_{X,Y}$, $c^{\Cgg}_{X,Y,Z}$, $c^{\Dgg}_{U,V,W}$, for all
$X,Y,Z\in\{X_0,\ldots,X_n\}$ and
$U,V,W\in\{\F(X_0),\ldots,\F(X_n)\}$, and possibly identity
functors. By an $\F$-iterate of multiplicity $n$,
or simply an $n$-iterate if there is no ambiguity, we will mean
a collection $H=\{ H_{X_0,\ldots,X_n}\}_{(X_0,\ldots,X_n)\in|\Cgg|^{n+1}}$,
where $H_{X_0,\ldots,X_n}$ is an $\F_{X_0,\ldots,X_n}$-iterate, called the
$(X_0,\ldots,X_n)$-component of $H$, the same for
all collections $X_0,\ldots,X_n$.
\end{defn}

\begin{rem}
{\rm When $\F$ is $K$-linear, the iterates may be thought of as $K$-linear
functors from $\Cgg(X_0,X_1)\odot\cdots\odot\Cgg(X_{n-1},X_n)$ to
$\Dgg(\F(X_0),\F(X_n))$, where $\odot$ denotes the Deligne product of
$K$-linear categories.}
\end{rem}

According to the previous definition, the image of the above path
$\gamma$ by the
$(X_0,\ldots,X_n)$-component
of a generic $n$-iterate $H$ will be of the form
$$
\F(f_n\circ\cdots\circ f_{i_1+\cdot+i_{r-1}+1})\circ\cdots\circ\F(f_{i_1+i_2+1}\circ\cdots\circ
f_{i_1+1})\circ\F(f_{i_1}\circ\cdots\circ f_1)
$$
for some ordered partition $\{1,\ldots,i_1\}$,
$\{i_1+1,\ldots,i_1+i_2\}$,$\{i_1+i_2+1,\ldots,i_1+i_2+i_3\}$,$\ldots$,
$\{i_1+\cdots+i_{r-1}+1,\ldots,i_1+i_2+\cdots+i_r\}$ of the set
$\{1,\ldots,n\}$, with $i_1+i_2+\cdots+i_r=n$ and $1\leq r\leq n$. Since such
a partition
completely defines the corresponding $n$-iterate and the partition itself
is completely given by the sequence $(i_1,\ldots,i_r)$, the
corresponding $n$-iterate will be denoted by $\F^{(i_1,\ldots,i_r)}$. For
example, there is a unique
$\F$-iterate of multiplicity $n=1$, namely, $\F^{(1)}$, given by the family of
functors defining the pseudofunctor $\F$ itself. For $n=2$, we have two
different $\F$-iterates, $\F^{(1,1)}$ and
$\F^{(2)}$, sending the path
$X\stackrel{f}{\To}Y\stackrel{g}{\To} Z$ to
$\F(g)\circ\F(f)$ and $\F(g\circ f)$, respectively. Their $(X,Y,Z)$-components
are given by
\begin{align*}
\F^{(1,1)}_{X,Y,Z}&=c^{\Dgg}_{\F(X),\F(Y),\F(Z)}\circ(\F_{X,Y}\times\F_{Y,Z})
\qquad X,Y\in|\Cgg| \\
\F^{(2)}_{X,Y,Z}&=\F_{X,Z}\circ c^{\Dgg}_{\F(X),\F(Y),\F(Z)} \qquad
X,Y\in|\Cgg|
\end{align*}
(in the $K$-linear case, the product $\times$ should be replaced by
the Deligne product $\odot$). The reader may easily check that there
are four 3-iterates,
which are exactly those defined by the families of functors appearing in
Lemma~\ref{axioma_composicio_reescrit}.

\begin{defn}
Given two $n$-iterates $H$,$H'$ of $\F$, $n\geq 1$, we will call {\sl indexed
  natural transformation} from $H$ to $H'$, and denote it by $\psi:H\Rightarrow H'$, any collection of natural
transformations between the corresponding components, i.e.,
$$
\psi=\{\psi_{X_0,\ldots,X_n}:H_{X_0,\ldots,X_n}\Rightarrow
H'_{X_0,\ldots,X_n}\}_{(X_0,\ldots,X_n)\in|\Cgg|^{n+1}}
$$
The natural transformation $\psi_{X_0,\ldots,X_n}$ will be called the
$(X_0,\ldots,X_n)$-component of $\psi$.
\end{defn}
Notice that, in this definition, no relation is required between the natural
transformations corresponding to the various components
$\psi_{X_0,\ldots,X_n}$ of $\psi$, for different collections of objects $(X_0,\ldots,X_n)$.

Given two such indexed
natural transformations $\psi:H\Rightarrow H'$ and $\psi':H'\Rightarrow H''$,
for some $n$-iterates $H,H',H''$, we define their vertical composite
as the indexed natural transformation $\psi'\cdot\psi:H\Rightarrow H''$ whose
components are given by the usual vertical composition of natural
transformations, i.e.,
\begin{equation} \label{composicio_vertical_trans_nat_indexades}
(\psi'\cdot\psi)_{X_0,\ldots,X_n}=\psi'_{X_0,\ldots,X_n}\cdot\psi_{X_0,\ldots,X_n}
\end{equation}

The 2-cosemisimplicial object of $\F$ in {\bf Cat} is then defined
as follows. Take $\Cc^0=\Cc^0(\F)={\bf 1}$, the terminal category
with only one object and one (identity) morphism. For $n\geq 1$,
let $\Cc^n(\F)$ be the small category with objects all $n$-iterates
of $\F$ and morphisms the indexed natural transformations between
them as defined above, the composition being the above vertical
composition. As regards the coface functors, they will be denoted
by $O_n^i:\Cc^{n-1}(\F)\rightarrow\Cc^n(\F)$, and they are defined
as follows. If $n=1$, both $O_1^0$ and $O_1^1$ are equal to the
unique possible functor from {\bf 1} to $\Cc^1(\F)$. If $n\geq 2$,
let

\begin{itemize}
\item
$O^0_n$ be the functor sending the $(n-1)$-iterate $H$ to
$$
O_n^0(H)_{X_0,\ldots,X_n}=c^{\Dgg}_{\F(X_0),\F(X_1),\F(X_n)}\circ
\left(\F_{X_0,X_1}\times H_{X_1,\ldots,X_n}\right)
$$
and an indexed natural transformation $\psi:H\Rightarrow H'$ to
$$
O_n^0(\psi)_{X_0,\ldots,X_n}={\bf 1}_{c^{\Dgg}_{\F(X_0),\F(X_1),\F(X_n)}}\circ
\left({\bf 1}_{\F_{X_0,X_1}}\times \psi_{X_1,\ldots,X_n}\right)
$$

\item
$O_n^i$, for $i=1,\ldots,n-1$, be the functor sending the $(n-1)$-iterate $H$ to
$$
O^i_n(H)_{X_0,\ldots,X_n}=H_{X_0,\ldots,\hat{X}_i,\ldots,X_n}\circ\left(\id_0
\times\cdots\times c^{\Cgg}_{X_{i-1},X_i,X_{i+1}}\times\cdots\times
\id_n\right)
$$
and an indexed natural transformation $\psi:H\Rightarrow H'$ to
$$
O^i_n(\psi)_{X_0,\ldots,X_n}=\psi_{X_0,\ldots,\hat{X}_i,\ldots,X_n}\circ\left({\bf
    1}_{\id_0}
\times\cdots\times {\bf 1}_{c^{\Cgg}_{X_{i-1},X_i,X_{i+1}}}\times\cdots\times
{\bf 1}_{\id_n}\right)
$$
(for short, we write here $\id_i$ instead of $\id_{\Cgg(X_i,X_{i+1})}$), and

\item
$O_n^n$ be the functor sending the $(n-1)$-iterate $H$ to
$$
O^n_n(H)_{X_0,\ldots,X_n}=
c^{\Dgg}_{\F(X_0),\F(X_{n-1}),\F(X_n)}\circ
\left(H_{X_0,\ldots,X_{n-1}}\times\F_{X_{n-1},X_n}\right)
$$
and an indexed natural transformation $\psi:H\Rightarrow H'$ to
$$
O^n_n(\psi)_{X_0,\ldots,X_n}=
{\bf 1}_{c^{\Dgg}_{\F(X_0),\F(X_{n-1}),\F(X_n)}}\circ
\left(\psi_{X_0,\ldots,X_{n-1}}\times{\bf 1}_{\F_{X_{n-1},X_n}}\right)
$$
\end{itemize}

The reader may easily check that the above formulas are
indeed functorial. Notice also that all these coface functors correspond to all
possible ways of
getting an $n$-iterate from an $(n-1)$-iterate.

It is a tedious but straightforward computation to check that these functors
$O^i_n$ satisfy the
cosemisimplicial identities {\rm (\ref{identitats_cosemisimplicials})}
for all $0\leq i<j\leq n+1$ except for the pairs $i=0,j=1$ and $i=n,j=n+1$,
with $n\geq 1$. When $n=1$, $O_2^1\circ
O_1^0:\Cc^0(\F)\rightarrow\Cc^2(\F)$ is the functor sending the unique object
$\star$ of $\Cc^0(\F)$ to
the 2-iterate $\F^{(2)}$, while $O_2^0\circ
O_1^0$ sends it to $\F^{(1,1)}$. Hence, it makes sense to define a natural
isomorphism $\tau^1_{0,1}:O_2^1\circ O_1^0\Rightarrow O_2^0\circ O_1^0$ whose
unique component $\tau^1_{0,1}(\star)$ is the indexed natural
transformation with $(X,Y,Z)$-component given by
$$
\tau^1_{0,1}(\star)_{X,Y,Z}:=\widehat{\F}^{-1}_{X,Y,Z}
$$
Similarly, $O_2^2\circ O_1^1$ sends the object $\star$ to the 2-iterate
$\F^{(1,1)}$ while $O_2^1\circ O_1^1$ sends it to $\F^{(2)}$, so that we can
define $\tau^1_{1,2}:O_2^2\circ O_1^1\Rightarrow O_2^1\circ O_1^1$ by
$$
\tau^1_{1,2}(\star)_{X,Y,Z}:=\widehat{\F}_{X,Y,Z}
$$
for all $X,Y,Z$. When $n\geq 2$, the images of an arbitrary
$(n-1)$-iterate $H$
by the functors $O^1_{n+1}\circ O^0_n$, $O^0_{n+1}\circ O^0_n$,
$O^{n+1}_{n+1}\circ O^n_n$ and $O^n_{n+1}\circ O^n_n$ are respectively given
by
\begin{align*}
(O^1_{n+1}\circ
  O^0_n)(H)_{X_0,\ldots,X_{n+1}}&=c^{\Dgg}_{\F(X_0),\F(X_2),\F(X_{n+1})}\circ\left(\F^{(2)}_{X_0,X_1,X_2}\times
  H_{X_2,\ldots,X_{n+1}}\right) \\
(O^0_{n+1}\circ
  O^0_n)(H)_{X_0,\ldots,X_{n+1}}&=c^{\Dgg}_{\F(X_0),\F(X_2),\F(X_{n+1})}\circ
\left(\F^{(1,1)}_{X_0X_1,X_2}
\times H_{X_2,\ldots,X_{n+1}}\right) \\
(O^{n+1}_{n+1}\circ
  O^n_n)(H)_{X_0,\ldots,X_{n+1}}&=c^{\Dgg}_{\F(X_0),\F(X_{n-1}),\F(X_{n+1})}\circ\left(H_{X_0,\ldots,X_{n-1}}\times\F^{(1,1)}_{X_{n-1},X_n,X_{n+1}}\right)
  \\ (O^n_{n+1}\circ
  O^n_n)(H)_{X_0,\ldots,X_{n+1}}&=c^{\Dgg}_{\F(X_0),\F(X_{n-1}),\F(X_{n+1})}\circ\left(H_{X_0,\ldots,X_{n-1}}\times\F^{(2)}_{X_{n-1},X_n,X_{n+1}}\right)
\end{align*}
Hence, for all $n\geq 2$, we can define natural isomorphisms
$\tau^n_{0,1}:O^1_{n+1}\circ O^0_n\Rightarrow O^0_{n+1}\circ O^0_n$ and
$\tau^n_{n,n+1}:O^{n+1}_{n+1}\circ O^n_n\Rightarrow O^n_{n+1}\circ O^n_n$
whose $H$-components, for any $(n-1)$-iterate
$H$, are the indexed natural transformations with components
$$
\tau^n_{0,1}(H)_{X_0,\ldots,X_{n+1}}={\bf
  1}_{c^{\Dgg}_{\F(X_0),\F(X_2),\F(X_{n+1})}}\circ
(\widehat{\F}^{-1}_{X_0,X_1,X_2}\times
  {\bf 1}_{H_{X_2,\ldots,X_{n+1}}})
$$
and
$$
\tau^n_{n,n+1}(H)_{X_0,\ldots,X_{n+1}}={\bf
  1}_{c^{\Dgg}_{\F(X_0),\F(X_{n-1}),\F(X_{n+1})}}\circ
({\bf 1}_{H_{X_0,\ldots,X_{n-1}}}\times\widehat{\F}_{X_{n-1},X_n,X_{n+1}})
$$
We have then the following:
\begin{thm}
For any pseudofunctor $\F:\Cgg\rightarrow\Dgg$, the triple
$(\Cc^{\bullet}(\F), O,\tau)$, with all $\tau^n_{ij}$'s equal to
identities except in the cases $(i=0,j=1)$ and $(i=n,j=n+1)$, where
they are given as above, defines a 2-cosemisimplicial object in
{\bf Cat}.
\end{thm}
\begin{proof}
We have to see that the 2-isomorphisms $\tau^n_{ij}$, as defined above,
satisfy the coherence relations in
Definition~\ref{objecte_2-cosemisimplicial}, for all triples $(i,j,k)$ with
$0\leq i<j<k\leq n+2$. Almost all such conditions are empty because
many of the $\tau$'s are trivial. It is easy to see that the only
nonempty conditions correspond to the triples $(i,j,k)$ of one of the
the following two families:
\begin{itemize}
\item
$i=0$, $j=1$ and $k\in\{2,\ldots,n+2\}$, and
\item
$i\in\{0,\ldots,n\}$, $j=n+1$ and $k=n+2$.
\end{itemize}
Let us consider the case $n=1$. In this case, the following four conditions
must be checked:
\begin{align*}
(1_{O_3^0}\circ\tau^1_{0,1})\cdot(\tau^2_{0,2}\circ
1_{O_1^0})\cdot(1_{O_3^2}\circ\tau^1_{0,1})&=(\tau^2_{0,1}\circ
1_{O_1^0})\cdot(1_{O_3^1}\circ\tau^1_{0,1})\cdot(\tau^2_{1,2}
\circ 1_{O_1^0}) \\ (1_{O_3^0}\circ\tau^1_{1,2})\cdot(\tau^2_{0,3}\circ
1_{O_1^1})\cdot(1_{O_3^3}\circ\tau^1_{0,2})&=(\tau^2_{0,2}\circ
1_{O_1^1})\cdot(1_{O_3^2}\circ\tau^1_{0,2})\cdot(\tau^2_{2,3}
\circ 1_{O_1^0}) \\ (1_{O_3^0}\circ\tau^1_{0,2})\cdot(\tau^2_{0,3}\circ
1_{O_1^0})\cdot(1_{O_3^3}\circ\tau^1_{0,1})&=(\tau^2_{0,1}\circ
1_{O_1^1})\cdot(1_{O_3^1}\circ\tau^1_{0,2})\cdot(\tau^2_{1,3}
\circ 1_{O_1^0}) \\ (1_{O_3^1}\circ\tau^1_{1,2})\cdot(\tau^2_{1,3}\circ
1_{O_1^1})\cdot(1_{O_3^3}\circ\tau^1_{1,2})&=(\tau^2_{1,2}\circ
1_{O_1^1})\cdot(1_{O_3^2}\circ\tau^1_{1,2})\cdot(\tau^2_{2,3}
\circ 1_{O_1^1})
\end{align*}
Proving any one of these equalities means checking that the
$\star$-component of both natural transformations (which are some
indexed natural transformation) coincide. The reader may easily check
that in the first and last cases, the condition one gets is the same,
namely
$$
\sigma^{24}_{X,Y,Z,T}\cdot\sigma^{12}_{X,Y,Z,T}=\sigma^{34}_{X,Y,Z,T}\cdot
\sigma^{13}_{X,Y,Z,T}
$$
where the natural transformations $\sigma^{ij}_{X,Y,Z,T}$ are those
defined in Lemma~\ref{axioma_composicio_reescrit}. Hence, both
conditions are equivalent to the composition axiom on $\F$. As regards
the other two equalities, they turn out to be true for all values of
$\widehat{\F}_{X,Y,Z,T}$. Indeed, the reader may check that
the $\star$-component of the left- and right-hand side natural
transformations in the second condition are both the indexed
natural transformation with components
$$
{\bf 1}_{c^{\Dgg}_{\F(X),\F(Y),\F(T)}}\circ\left({\bf
  1}_{\F_{X,Y}}\times\widehat{\F}_{Y,Z,T}\right)
$$
while in the third condition both are the indexed natural
transformation defined by
$$
{\bf
  1}_{c^{\Dgg}_{\F(X),\F(Z),\F(T)}}\circ\left(\widehat{\F}^{-1}_{X,Y,Z}\times {\bf 1}_{\F_{Z,T}}\right)
$$
When $n\geq 2$, the situation is similar. For the extreme values
$k=2,n+2$ it
turns out that both conditions reduce to the composition axiom on
$\F$, while in the cases $k=3,\ldots,n+1$ they are always satisfied,
for all values of $\widehat{\F}_{X,Y,Z,T}$.
\end{proof}

Suppose now that $\F$ is $K$-linear. The $K$-linear structure
in the target 2-category $\Dgg$ naturally induces a $K$-module
structure on the set of indexed natural transformations between any
two iterates so that the categories $\Cc^n(\F)$ are
$K$-linear. Furthermore, from
the definition of Deligne product of natural transformations between
$K$-linear functors it immediately follows that all coface
functors $O_n^i$ are also $K$-linear. Hence, the
corresponding 2-cosemisimplicial object of $\F$ belongs in this case to ${\bf
  Cat}_K$. Furthermore, it is trivially enhanced, because
$O_1^0=O_1^1$, so that the cochain complex construction of the previous section
can be applied. Notice that, in this case, we have no choice for the object
$X\in|\Cc^0|$, because $\Cc^0(\F)$ has only one object.

\begin{prop}
If $\F$ is $K$-linear, its deformation complex $X^{\bullet}(\F)$ coincides
 with the cochain complex obtained from the previous 2-cosemisimplicial
 object by the method described above when we take as reference objects in
 $\Cc^n(\F)$, for $n\geq 1$, those defined
inductively by Equations
(\ref{objectes_referencia_1})-(\ref{objectes_referencia_2}).
\end{prop}
\begin{proof}
It is easy
to see that these reference objects indeed correspond to the $n$-iterates used
in Section 3 to define $X^{\bullet}(\F)$, i.e.
\begin{align*}
X_n&=\F^{(1,\stackrel{n)}{\ldots},1)}
\\ X'_n&=\F^{(n)}
\end{align*}
The corresponding $K$-modules $M^n={\rm Hom}_{\Cc^n(\F)}
(\F^{(1,\stackrel{n)}{\ldots},1)},\F^{(n)})$ may then be identified
with the $X^n(\F)$ defined in Section 3. Moreover, under this identification,
the coboundary maps given by Equation
(\ref{operadors_covora}) exactly correspond to those defined in
Section 3 for the $K$-modules
$X^{\bullet}(\F)$, the action of the padding operators corresponding
to taking the left and right composites with the canonical isomorphisms
$\alpha_{i,n}^{{\bf 1}}$ and $\beta_{i,n}^{{\bf 1}}$.
\end{proof}


\section{Deviation calculus for an arbitrary $K$-linear category}

\label{seccio_desviacio}

A basic question regarding the deformation theory of a $K$-linear
pseudofunctor which remained open in \cite{jE1} is that of the
higher-order obstructions. To settle down this question, we
introduce in this section a generalization to arbitrary $K$-linear
categories of the deviation calculus introduced by Markl and
Stasheff for the category of $K$-modules \cite{MS94} and state the
corresponding additivity principle. The 2-cosemisimplicial object
of $\F$ introduced in the previous section turns out to fit quite
naturally in the framework of this deviation calculus and allows us
to give an easy proof that the higher-order obstructions are indeed
cocycles in the deformation complex. The proof is deferred to the
next section.

Let us start with the following definition, which generalizes the
$K[[h]]$-linear extension of a $K$-linear category and provides the right
setting in which doing a deviation calculus.
\begin{defn}
Let $\Cc$ be any $K$-linear category. Then, we will call {\sl deviation
  extension} of $\Cc$ any complete $K[[h]]$-linear category $\Cc_h$
  which is $K[[h]]$-linear isomorphic
  to the $K[[h]]$-linear extension $\Cc[[h]]$ defined above.
\end{defn}

Hence, if $\Cc_h$ is a deviation category of $\Cc$, we have a
bijection between objects
$\varphi:|\Cc_h|\rightarrow|\Cc|$ and $K[[h]]$-linear isomorphisms
$\Cc_h(X_h,Y_h)\cong\Cc(X,Y)[[h]]$ for all $X_h,Y_h\in|\Cc_h|$
(where $X=\varphi(X_h)$ and $Y=\varphi(Y_h)$) such that the
composition of morphisms
in $\Cc_h$ corresponds, after these identifications, to taking the usual
``product'' of formal power series.

\begin{ex}
{\rm If $\Cc=Mod_K$, the category of $K$-modules, then the full subcategory
$Mod_{K[[h]]}^0$ of $Mod_{K[[h]]}$ with objects the topologically free $K[[h]]$-modules is a deviation
extension of $\Cc$. This follows from the well-known isomorphisms of
$K[[h]]$-modules ${\rm
  Hom}_{K[[h]]}(V[[h]],W[[h]])\cong ({\rm Hom}_K(V,W))[[h]]$.}
\end{ex}
This is the example considered by Markl and Stasheff. The example we
are interested in this paper is the following.
\begin{ex} \label{exemple_categoria_desviacio}
{\rm  Let
$\Aa$, $\Bb$ be
two $K$-linear categories. Then, the functor
category ${\rm Fun}_K(\Aa,\Bb)$ with objects all $K$-linear functors
$F:\Aa\rightarrow\Bb$ and morphisms the natural transformations with the
vertical composition is a $K$-linear category. It turns out that a
deviation extension of ${\rm Fun}_K(\Aa,\Bb)$ is given by the full
subcategory of the functor category ${\rm
  Fun}_{K[[h]]}(\Aa[[h]],\Bb[[h]])$ with objects all $K[[h]]$-linear
functors $F_h:\Aa[[h]]\rightarrow\Bb[[h]]$ which are $K[[h]]$-linear
extensions $F_h=F[[h]]$ of a $K$-linear functor
$F:\Aa\rightarrow\Bb$. Let us denote this subcategory by ${\rm
  Fun}_{K[[h]]}(\Aa[[h]],\Bb[[h]])^0$. That such a
category is a deviation extension of ${\rm Fun}_K(\Aa,\Bb)$ follows from
Lemma~\ref{lema_nat_h}. }
\end{ex}

Let $\Cc_h$ be a deviation category of $\Cc$ and let us fix isomorphisms
$\varphi_{X,Y}$ as above. We will identify each morphism in $\Cc_h$ with
the corresponding formal power series as given by these isomorphisms. Let us
then consider a ``potentially commutative'' diagram in $\Cc_h$ of the form
$$
\xymatrix{
X_h\ar[r]^{\alpha_h}
\ar[d]_{\gamma_h} &
Y_h\ar[d]^{\beta_h} \\ T_h\ar[r]_{\delta_h} & Z_h }
$$
with $\alpha_h=\sum_{n\geq 0}\alpha_n h^n$, $\alpha_n\in\Cc(X,Y)$ for all
$n\geq 0$, and similarly $\beta_h$, $\gamma_h$
and $\delta_h$. Since the composition of two consecutive morphisms in this
diagram is given by
the usual product rule between formal power series, the commutativity of the
diagram is equivalent to the infinite set of equations
$$
\sum_{p+q=m}(\beta_p\circ\alpha_q-\delta_p\circ\gamma_q)=0,\quad m\geq 0.
$$
Hence, it makes sense to talk about the commutativity of such a
diagram modulo $h^{n+1}$ (the equations are satisfied for all $m\leq
n$ but possibly not for $m=n+1$). Following Markl and Stasheff
\cite{MS94}, one may
then define the deviation for such a diagram as the first non
zero coefficient of the map $\delta_h\circ\gamma_h-\beta_h\circ\alpha_h$. More
explicitly:

\begin{defn}
Suppose that a potentially commutative diagram in $\Cc_h$ as above
commutes modulo $h^{n+1}$, but not modulo $h^{n+2}$. Then, the {\sl
  deviation} of the diagram is the
(unique) morphism $\Psi:X\To Z$ in $\Cc$ determined by the equation
$$
\delta_h\circ\gamma_h-\beta_h\circ\alpha_h=\Psi h^{n+1} \quad {\rm mod}\
h^{n+2}
$$
\end{defn}

\begin{rem}
{\rm A priori, the deviation as defined here may depend on the isomorphisms
$\varphi_{X,Y}$ giving $\Cc_h$ the structure of a deviation extension of
$\Cc$. This is the reason by which we need to fix these isomorphisms.}
\end{rem}

\begin{ex}
{\rm Given a $K$-linear pseudofunctor $\F$, let $\Cc={\rm
  Fun}_K(\Cc^2(\F),\Cc^3(\F))$, where $\Cc^2(\F)$ and $\Cc^3(\F)$ are the
  categories that appear in the definition of the 2-cosemisimplicial
  object associated to $\F$. This is a $K$-linear
  category of the form considered in
  Example~\ref{exemple_categoria_desviacio} and a diagram in the corresponding
  deviation extension is precisely the collection of diagrams
(\ref{diagrama_desviacio=obstruccio}) appearing in
  Lemma~\ref{condicions_deformacio}, for all objects $X,Y,Z,T$. If such
  diagrams commute modulo $h^{n+1}$ but not modulo $h^{n+2}$, an easy
  computation gives that the deviation is the indexed natural transformation
 $\Psi$ with components
\begin{equation} \label{obstruccio}
\begin{array}{l}
\Psi_{X,Y,Z,T}=\sum_{p+q=n+1}
\left[\widehat{\F}^p_{X,Z,T}\circ{\mathbf 1}_{c^{\Cgg}_{X,Y,Z}\odot
\id_{\Cgg(Z,T)}}\right]\cdot \\ \hspace{4 truecm}\cdot\left[{\mathbf
  1}_{c^{\Dgg}_{\F(X),\F(Z),\F(T)}}\circ\left(\widehat{\F}^q_{X,Y,Z}\odot{\mathbf 1}_{\F_{Z,T}}\right)\right] \\ \\
\hspace{1.7 truecm}-\sum_{p+q=n+1}
\left[\widehat{\F}^p_{X,Y,T}\circ{\mathbf 1}_{\id_{\Cgg(X,Y)}\odot
c^{\Cgg}_{Y,Z,T}}\right]\cdot \\ \hspace{4 truecm}\cdot\left[{\mathbf 1}_{c^{\Dgg}_{\F(X),\F(Y),\F(T)}}\circ\left({\mathbf 1}_{\F_{X,Y}}\odot\widehat{\F}^q_{Y,Z,T}\right)\right]
\end{array}
\end{equation} }
\end{ex}

Notice that, in the previous definition, one implicitly chooses an order
between the two paths in the
diagram, and that the same diagram with the reverse order
corresponds to the same deviation but with opposite sign. To
indicate which deviation one is considering, an arrow is sometimes drawn in
the diagram from the first to the second path. In the example above, $\Psi$ is
the deviation from the path $\sigma^{24}(h)\cdot\sigma^{12}(h)$ to the path
$\sigma^{34}(h)\cdot\sigma^{13}(h)$. Clearly, the definition may be extended
without trouble to the deviation of
any potentially commutative diagram of an arbitrary polygonal shape.

The fundamental point in Markl and Stasheff's deviation calculus is
an easy {\it additivity principle} which allows one to compute the deviation
of any potentially commutative diagram having the form of a polygonally
subdivided diagram such as that below.

$$
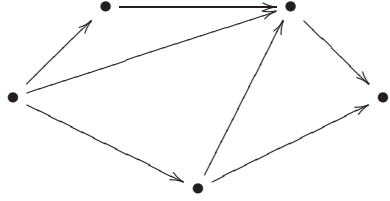
\xymatrix{
& \bullet\ar[rr]^{} & & \bullet\ar[rd]^{}& \\
\bullet\ar[ru]^{}\ar[rrru]^{}
\ar[rrd]_{} & & & & \bullet \\
& & \bullet\ar[ruu]_{}\ar[rru]_{} & &  }
$$

\noindent{In} our general context, this principle can be stated as follows:

\begin{prop} \label{principi_additivitat}
Let $\Cc_h$ be a deviation category of $\Cc$ with fixed isomorphisms
$\varphi_{X,Y}$ for all $X,Y\in|\Cc|$, and let us consider two diagrams
in $\Cc_h$ with a common edge
$$
\begin{array}{ccc}
\xymatrix{
X_h\ar[r]^{\alpha_h}
\ar[d]_{\gamma_h} &
Y_h\ar[d]^{\beta_h} \\ T_h\ar[r]_{\delta_h} & Z_h } & &
\xymatrix{
Y_h\ar[r]^{\epsilon_h}
\ar[d]_{\beta_h} &
U_h\ar[d]^{\eta_h} \\ Z_h\ar[r]_{\xi_h} & V_h }
\end{array}
$$
Suppose that both diagrams commute modulo $h^{n+1}$ but not modulo
$h^{n+2}$ and let $\Psi_1:X\To Z$, $\Psi_2:Y\To V$ denote the corresponding
deviations. Then, the composite diagram
$$\xymatrix{
X_h\ar[r]^{\epsilon_h\circ\alpha_h}
\ar[d]_{\gamma_h} &
U_h\ar[d]^{\eta_h} \\ T_h\ar[r]_{\xi_h\circ\delta_h} & V_h }
$$
commutes modulo $h^{n+1}$ but not modulo $h^{n+2}$ and its deviation
$\Psi:X\To V$ is given by
$$
\Psi=\xi_0\circ\Psi_1+\Psi_2\circ\alpha_0
$$
\end{prop}
\begin{proof}
The proof is formally the same as in the case $\Cc=Mod_K$ and is left to the
reader (see \cite{MS94}).
\end{proof}
Note that, when the zero order terms of the maps $\alpha_h$ and
$\xi_h$ are identities (in particular, $Y=X$ and $V=Z$), deviations simply
add, suggesting the name ``additivity principle'' for this result.
Using such a basic additivity principle, we can easily get expressions for
the deviation of more complex diagrams. For example, the reader may easily
check that the deviation of the previous diagram
is simply given by the sum of the deviations of each of the three faces.

For our purposes, the relevant result on this
deviation calculus is the following obvious fact:

\vspace{0.4 truecm}
\noindent{{\bf Basic fact.}}
{\it Let $D_1,D_2$ be two potentially commutative subdivided polygonal
diagrams in a deviation extension $\Cc_h$ of a $K$-linear category $\Cc$,
commuting modulo $h^{n+1}$ and
with a common boundary (consequently defining a
2-dimensional polyhedron topologically equivalent to $S^2$). Then,
the deviations of both diagrams must coincide.}


\section{Higher-order obstructions}

Let us now consider the question of the obstructions. Our purpose in this
section is to prove, using the previous deviation calculus, that the
obstruction to the integrability of a purely pseudofunctorial $n^{th}$-order
deformation of $\F$ indeed corresponds to a cocycle in the
deformation complex. More explicitly, we have the following.

\begin{thm}
The obstruction to the extension one higher order of a purely
pseudofunctorial $n^{th}$-order deformation of a $K$-linear
unitary pseudofunctor $\F$ is a 3-cocycle in the corresponding deformation
complex $X^{\bullet}(\F)$. If this obstruction cocycle defines the zero
cohomolofy class in $H^3(\F)$ an extension exists.
\end{thm}

\begin{proof}
Let
$(\widehat{\F}_h)_{X,Y,Z}=\widehat{\F}_{X,Y,Z}+\widehat{\F}^1_{X,Y,Z}h+\cdots+\widehat{\F}^n_{X,Y,Z}h^n$
be a purely pseudofunctorial $n^{th}$-order deformation of
$\F$. Given
$\widehat{\F}^{n+1}=\{\widehat{\F}^{n+1}_{X,Y,Z}\}_{X,Y,Z}\in X^2(\F)$,
  an easy degree computation shows that
$(\widehat{\F}_h)_{X,Y,Z}+\widehat{\F}^{n+1}_{X,Y,Z}h^{n+1}$ defines
an $(n+1)$-deformation of the same kind if and only if
\begin{equation} \label{equacio_obstruccio}
\delta(\widehat{\F}^{n+1})=\Psi
\end{equation}
where the obstruction $\Psi=\{\Psi_{X,Y,Z,T}:\F^{(1,1,1)}_{X,Y,Z,T}\Rightarrow
\F^{(3)}_{X,Y,Z,T}\}_{X,Y,Z,T}\ \in X^3(\F)$ is the indexed natural
transformation with components
\begin{equation} \label{obstruccio_correcta}
\begin{array}{l}
\Psi_{X,Y,Z,T}=\sum_{\begin{array}{c} p+q=n+1 \\ 0\leq p,q\leq n \end{array}}
\left[\widehat{\F}^p_{X,Z,T}\circ{\mathbf 1}_{c^{\Cgg}_{X,Y,Z}\odot
\id_{\Cgg(Z,T)}}\right]\cdot \\ \hspace{4 truecm}\cdot\left[{\mathbf
  1}_{c^{\Dgg}_{\F(X),\F(Z),\F(T)}}\circ\left(\widehat{\F}^q_{X,Y,Z}
\odot{\mathbf 1}_{\F_{Z,T}}\right)\right] \\ \\
\hspace{1.7 truecm}-\sum_{\begin{array}{c} p+q=n+1 \\ 0\leq p,q\leq n
\end{array}}
\left[\widehat{\F}^p_{X,Y,T}\circ{\mathbf 1}_{\id_{\Cgg(X,Y)}\odot
c^{\Cgg}_{Y,Z,T}}\right]\cdot \\ \hspace{4 truecm}\cdot
\left[{\mathbf 1}_{c^{\Dgg}_{\F(X),\F(Y),\F(T)}}\circ
\left({\mathbf 1}_{\F_{X,Y}}\odot\widehat{\F}^q_{Y,Z,T}\right)\right]
\end{array}
\end{equation}
Notice that these are exactly the components of the indexed natural
transformation giving the deviation of
diagrams (\ref{condicions_deformacio}) (see Equation~(\ref{obstruccio}))
except that the sums are taken over all $p+q=n+1$ such that $0\leq
p,q\leq n$. Such restrictions are due to the fact that we are now considering
the deviation of diagrams (\ref{condicions_deformacio}) when the
$\sigma^{ij}(h)$
are those defined by the $n^{th}$-order deformation $\widehat{\F}_h$ (in
particular, we indeed have commutativity modulo $h^{n+1}$).

We want to see that
$\delta(\Psi)=0$ (this is the necessary condition for an
$\widehat{\F}^{n+1}$ satisfying Equation~(\ref{equacio_obstruccio}) to
exist). From the definition of $\delta:X^2(\F)\rightarrow X^3(\F)$ as given in
Equation~(\ref{operadors_covora}), we have that
\begin{equation} \label{delta_psi}
\delta(\Psi)=\sum_{i=0}^4 (-1)^i\alpha_{i,4}\cdot O_4^i(\Psi)\cdot\beta_{i,4}
\end{equation}
where $\alpha_{i,4}$ and $\beta_{i,4}$, $i=0,1,2,3,4$, denote the
$\star$-components of some canonical enhanced natural
isomorphisms. Explicitly,
\begin{align*}
\alpha_{i,4}&=\left(\tau_{(0,1,2,i),(0,1,2,3)}\right)_{\star} \\
\beta_{i,4}&=\left(\tau_{(1,2,3,i),(1,2,3,4)}\right)_{\star}
\end{align*}
(although not made explicit in the $\alpha$'s and $\beta$'s, recall that we
are taking as enhancing isomorphism $\phi$ the identity
natural transformation of $O_1^0=O_1^1$). Notice that the composition in this
case is denoted by a dot because it corresponds to the vertical composition of
indexed natural transformations (see
Equation~(\ref{composicio_vertical_trans_nat_indexades})).

To prove that this is indeed the zero indexed natural transformation, let us
apply
the $K[[h]]$-linear extensions of the functors $O_4^0$, $O_4^1$, $O_4^2$,
$O_4^3$ and $O_4^4$ to the diagram~(\ref{axioma_hexagon}) from which $\Psi$
is the deviation \footnote{Strictly, what we apply to this diagram are
  not the $K[[h]]$-linear extensions of the $O_4^i$, because such extensions
  act on the category $\Cc^3(\F)[[h]]$, whose objects are the same as in
  $\Cc^3(\F)$. But we need to consider a category whose objects are the
  $K[[h]]$-linear extensions of the 3-iterates, not the 3-iterates
  themselves. Anyway, the meaning of these slightly different versions of the
  $O_n^i[[h]]$ is obvious.}. We leave to the reader to check that one obtains
the following new diagrams.

\vspace{0.3 truecm}
Action of $O_4^0[[h]]$:
\begin{equation} \label{O_0}
\xymatrix{
\F^{(1,1,1,1)}_{X,Y,Z,T}[[h]]\ar@{=>}[rr]^{\sigma^{1,4}_{X,Y,Z,T}(h)}
\ar@{=>}[d]_{\sigma^{2,4}_{X,Y,Z,T}(h)} & &
\F^{(1,2,1)}_{X,Y,Z,T}[[h]]\ar@{=>}[d]^{\sigma^{5,4}_{X,Y,Z,T}(h)} \\
\F^{(1,1,2)}_{X,Y,Z,T}[[h]]\ar@{=>}[rr]_{\sigma^{7,4}_{X,Y,Z,T}(h)} & & \F^{(1,3)}_{X,Y,Z,T}[[h]] }
\end{equation}

Action of $O^1_4[[h]]$:
\begin{equation} \label{O_1}
\xymatrix{
\F^{(2,1,1)}_{X,Y,Z,T}[[h]]\ar@{=>}[rr]^{\sigma^{4,4}_{X,Y,Z,T}(h)}
\ar@{=>}[d]_{\sigma^{3,4}_{X,Y,Z,T}(h)} & &
\F^{(3,1)}_{X,Y,Z,T}[[h]]\ar@{=>}[d]^{\sigma^{11,4}_{X,Y,Z,T}(h)} \\
\F^{(2,2)}_{X,Y,Z,T}[[h]]\ar@{=>}[rr]_{\sigma^{10,4}_{X,Y,Z,T}(h)} & & \F^{(4)}_{X,Y,Z,T}[[h]] }
\end{equation}

Action of $O^2_4[[h]]$:
\begin{equation} \label{O_2}
\xymatrix{
\F^{(1,2,1)}_{X,Y,Z,T}[[h]]\ar@{=>}[rr]^{\sigma^{6,4}_{X,Y,Z,T}(h)}
\ar@{=>}[d]_{\sigma^{5,4}_{X,Y,Z,T}(h)} & &
\F^{(3,1)}_{X,Y,Z,T}[[h]]\ar@{=>}[d]^{\sigma^{11,4}_{X,Y,Z,T}(h)} \\
\F^{(1,3)}_{X,Y,Z,T}[[h]]\ar@{=>}[rr]_{\sigma^{9,4}_{X,Y,Z,T}(h)} & & \F^{(4)}_{X,Y,Z,T}[[h]] }
\end{equation}

Action of $O^3_4[[h]]$:
\begin{equation} \label{O_3}
\xymatrix{
\F^{(1,1,2)}_{X,Y,Z,T}[[h]]\ar@{=>}[rr]^{\sigma^{8,4}_{X,Y,Z,T}(h)}
\ar@{=>}[d]_{\sigma^{7,4}_{X,Y,Z,T}(h)} & &
\F^{(2,2)}_{X,Y,Z,T}[[h]]\ar@{=>}[d]^{\sigma^{10,4}_{X,Y,Z,T}(h)} \\
\F^{(1,3)}_{X,Y,Z,T}[[h]]\ar@{=>}[rr]_{\sigma^{9,4}_{X,Y,Z,T}(h)} & & \F^{(4)}_{X,Y,Z,T}[[h]] }
\end{equation}

Action of $O^4_4[[h]]$:
\begin{equation} \label{O_4}
\xymatrix{
\F^{(1,1,1,1)}_{X,Y,Z,T}[[h]]\ar@{=>}[rr]^{\sigma^{0,4}_{X,Y,Z,T}(h)}
\ar@{=>}[d]_{\sigma^{1,4}_{X,Y,Z,T}(h)} & &
\F^{(2,1,1)}_{X,Y,Z,T}[[h]]\ar@{=>}[d]^{\sigma^{4,4}_{X,Y,Z,T}(h)} \\
\F^{(1,2,1)}_{X,Y,Z,T}[[h]]\ar@{=>}[rr]_{\sigma^{6,4}_{X,Y,Z,T}(h)} & & \F^{(3,1)}_{X,Y,Z,T}[[h]] }
\end{equation}
where $\F^{(1,1,1,1)}$, $\F^{(2,1,1)}$, $\F^{(1,2,1)}$, $\F^{(1,3)}$,
$\F^{(2,2,)}$, $\F^{(3,1)}$ and $\F^{(4)}$ denote the eight 4-iterates of $\F$
and the $\sigma^{i,4}_{X,Y,Z,T,U}(h)$ are formal power series in $h$ of the
form
$$
\sigma^{i,4}_{X,Y,Z,T,U}(h)=\sum_{k\geq 0} (\sigma^{i,4}_k)_{X,Y,Z,T,U} h^k
$$
with $\sigma^{i,4}_k=\{(\sigma^{i,4}_k)_{X,Y,Z,T,U}\}_{X,Y,Z,T,U}$, for $0\leq
i\leq 11$ and $k\geq 0$, the indexed natural transformations with components

\begin{align*}
(\sigma_k^{0,4})_{X,Y,Z,T,U}&={\mathbf
      1}_{c^{\Dgg}_{\F(X),\F(Z),\F(T),\F(U)}}\circ\left(\widehat{\F}^k_{X,Y,Z}\odot{\mathbf
      1}_{\F_{Z,T}}\odot{\mathbf 1}_{\F_{T,U}}\right)
\\
(\sigma^{1,4}_k)_{X,Y,Z,T,U}&={\mathbf
      1}_{c^{\Dgg}_{\F(X),\F(Y),\F(T),\F(U)}}\circ\left({\mathbf
      1}_{\F_{X,Y}}\odot\widehat{\F}^k_{X,Y,Z}\odot{\mathbf
      1}_{\F_{T,U}}\right)
\\
(\sigma^{2,4}_k)_{X,Y,Z,T,U}&={\mathbf
      1}_{c^{\Dgg}_{\F(X),\F(Y),\F(Z),\F(U)}}\circ\left({\mathbf
      1}_{\F_{X,Y}}\odot{\mathbf 1}_{\F_{Y,Z}}\odot
\widehat{\F}^k_{Z,T,U}\right)
\\
(\sigma^{3,4}_k)_{X,Y,Z,T,U}&={\mathbf
      1}_{c^{\Dgg}_{\F(X),\F(Z),\F(U)}}\circ\left({\mathbf
      1}_{\F_{X,Z}}\odot\widehat{\F}^k_{Z,T,U}\right)\circ{\mathbf 1}_{c^{\Cgg}_{X,Y,Z}\odot\id_{Z,T,U}}
\\
(\sigma^{4,4}_k)_{X,Y,Z,T,U}&={\mathbf
      1}_{c^{\Dgg}_{\F(X),\F(T),\F(U)}}\circ\left(\widehat{\F}^k_{X,Z,T}\odot{\mathbf
      1}_{\F_{T,U}}\right)\circ{\mathbf 1}_{c^{\Cgg}_{X,Y,Z}\odot\id_{Z,T,U}}
\\
(\sigma^{5,4}_k)_{X,Y,Z,T,U}&={\mathbf
      1}_{c^{\Dgg}_{\F(X),\F(Y),\F(U)}}\circ\left({\mathbf
      1}_{\F_{X,Y}}\odot\widehat{\F}^k_{Y,T,U}\right)\circ{\mathbf
      1}_{\id_{X,Y}\odot c^{\Cgg}_{Y,Z,T}\odot\id_{T,U}}
\\
(\sigma^{6,4}_h)_{X,Y,Z,T,U}&={\mathbf
      1}_{c^{\Dgg}_{\F(X),\F(T),\F(U)}}\circ\left(\widehat{\F}^k_{X,Y,T}\odot{\mathbf
      1}_{\F_{T,U}}\right)\circ{\mathbf
      1}_{\id_{X,Y}\odot
      c^{\Cgg}_{Y,Z,T}\odot\id_{T,U}}
\\
(\sigma^{7,4}_h)_{X,Y,Z,T,U}&={\mathbf
      1}_{c^{\Dgg}_{\F(X),\F(Y),\F(U)}}\circ\left({\mathbf
      1}_{\F_{X,Y}}\odot\widehat{\F}^k_{Y,Z,U}\right)\circ{\mathbf
      1}_{\id_{X,Y,Z}\odot c^{\Cgg}_{Y,Z,T}}
\\
(\sigma^{8,4}_k)_{X,Y,Z,T,U}&={\mathbf
      1}_{c^{\Dgg}_{\F(X),\F(Z),\F(U)}}\circ\left(\widehat{\F}^k_{X,Y,Z}\odot
{\mathbf 1}_{\F_{Z,U}}\right)\circ{\mathbf
      1}_{\id_{X,Y,Z}\odot c^{\Cgg}_{Y,Z,T}}
\\
(\sigma^{9,4}_k)_{X,Y,Z,T,U}&=\widehat{\F}^k_{X,Y,U}\circ{\mathbf
      1}_{\id_{\Cgg(X,Y)}\odot c^{\Cgg}_{Y,Z,T,U}}
\\
(\sigma^{10,4}_k)_{X,Y,Z,T,U}&=\widehat{\F}^k_{X,Z,U}\circ{\mathbf
      1}_{c^{\Cgg}_{X,Y,Z}\odot c^{\Cgg}_{Z,T,U}}
\\
(\sigma^{11,4}_k)_{X,Y,Z,T,U}&=\widehat{\F}^k_{X,T,U}\circ{\mathbf
      1}_{ c^{\Cgg}_{X,Y,Z,T}\odot\id_{\Cgg(T,U)}}
\end{align*}
%
%

These diagrams are five of the six faces of the cube in
Fig.~\ref{cub_pseudofunctor} (for short, when naming the vertices and edges in
this diagram, the indexing objects and the formal parameter $h$ have been
omitted).

\begin{figure}[hp]
\centering
\input{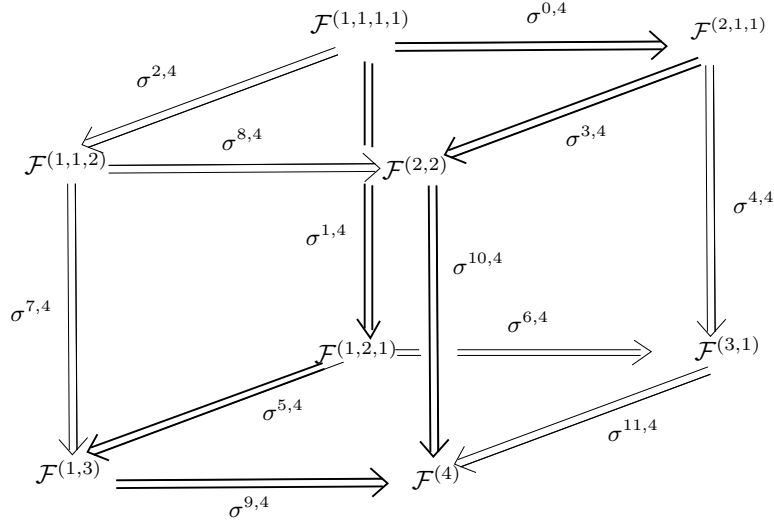}
\caption{Action of the functors $O^i_4$, $i=0,1,2,3,4$ on the
  diagram~(\ref{axioma_hexagon}).}
\label{cub_pseudofunctor}
\end{figure}

\begin{figure}[hp]
\centering
\input{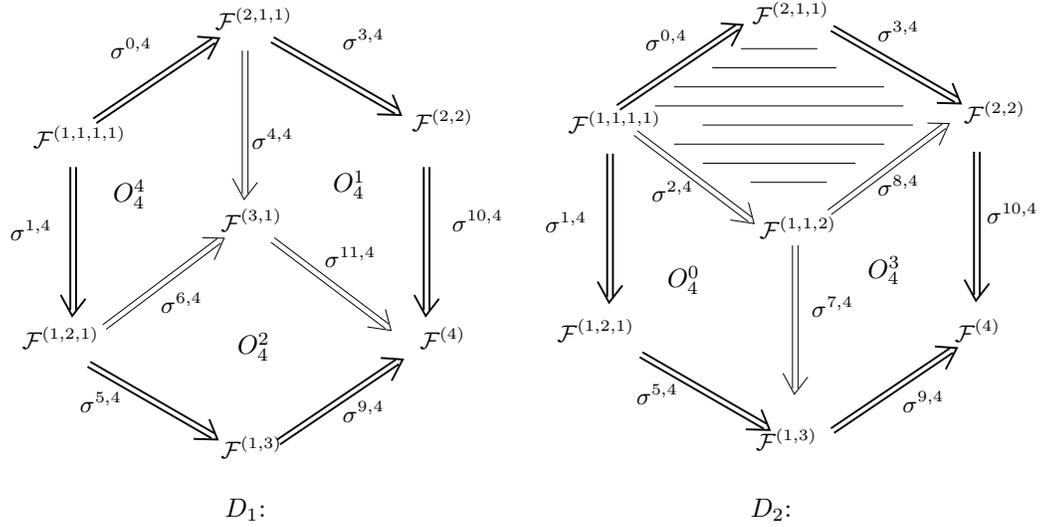}
\caption{Diagrams $D_1$ and $D_2$ decomposing the cube in Fig. 5.}
\label{D1_D2}
\end{figure}

As regards the lacking face at the top, it turns out to be always
commutative
(hence, it has null deviation). Indeed, the reader may easily check that, for
any path of 1-morphisms
$X\stackrel{f}{\longrightarrow}Y\stackrel{g}{\longrightarrow}
Z\stackrel{l}{\longrightarrow} T\stackrel{m}{\longrightarrow} U$ in $\Cgg$,
the $(f,g,l,m)$-component of the
degree $n$ term in the formal power series giving the composite
$\sigma^{3,4}_{X,Y,Z,T,U}(h)\cdot\sigma^{0,4}_{X,Y,Z,T,U}(h)$ is the
2-morphism
$$
\sum_{p+q=n} \left(\widehat{\F}^p(l,m)\circ 1_{\F(g\circ
  f)}\right)\cdot\left( 1_{\F(m)\circ\F(l)}\circ\widehat{\F}^q(f,g)\right)
$$
while the same component of the same term for
$\sigma^{8,4}_{X,Y,Z,T,U}(h)\cdot\sigma^{2,4}_{X,Y,Z,T,U}(h)$ is
$$
\sum_{p+q=n} \left( 1_{\F(m\circ l)}\circ\widehat{\F}^q(f,g)\right)\cdot\left(\widehat{\F}^p(l,m)\circ 1_{\F(g)\circ\F(f)}\right)
$$
By the interchange law, however, both 2-morphisms coincide  with
$\widehat{\F}^p(l,m)\circ\widehat{\F}^q(f,g)$, so that both composites are
equal. Hence, the above diagrams nicely fit in a 3-dimensional diagram {\bf D}
topologically
equivalent to $S^2$ and to which the basic fact from
Section~\ref{seccio_desviacio} may be applied. Looking at this diagram, it can
clearly be subdivided into the two hexagonal diagrams $D_1$ and $D_2$ depicted
in Fig.~\ref{D1_D2}, whose common boundary is indicated by bold arrows in
Fig.~\ref{cub_pseudofunctor}. Using now the additivity principle
(Proposition~\ref{principi_additivitat}), one obtains for the deviation of
$D_1$ the indexed natural transformation
$$
{\rm Dev}(D_1)=\sigma^{11,4}_0\cdot O_4^4(\Psi)-O_1^4(\Psi)\cdot\sigma^{0,4}_0+O_2^4(\Psi)\cdot\sigma^{1,4}_0
$$
while the deviation of $D_2$ turns out to be
$$
{\rm Dev}(D_2)=-\sigma^{9,4}_0\cdot
O_0^4(\Psi)+O_3^4(\Psi)\cdot\sigma^{2,4}_0
$$
By the basic fact in the previous section, we know that
$$
{\rm Dev}(D_1)={\rm Dev}(D_2)
$$
We leave to the reader to check that this is exactly the condition
$\delta(\Psi)=0$. Notice that taking the composites with the terms
  $\sigma^{i,4}_0$ in the above expressions for the deviations of $D_1$ and
  $D_2$, as established in the additivity principle, corresponds to taking the
  composites with the
  $\alpha_{i,4}$'s and $\beta_{i,4}$'s in Equation~(\ref{delta_psi}) and hence
  to the action of the padding operators.
\end{proof}

\vspace{0.7 truecm}
\noindent{{\bf Acknowledgements.}} I would like to thank L. Crane and
G. Rodrigues for calling my attention to R. Street's works on
descent categories. I am also indebted to the referee for his very
illuminating comments concerning Section 4.


\bibliographystyle{amsplain}
\bibliography{def_pseudof_def}

\providecommand{\bysame}{\leavevmode\hbox to3em{\hrulefill}\thinspace}
\providecommand{\MR}{\relax\ifhmode\unskip\space\fi MR }
\providecommand{\MRhref}[2]{%
  \href{http://www.ams.org/mathscinet-getitem?mr=#1}{#2}
}
\providecommand{\href}[2]{#2}
\begin{thebibliography}{10}

\bibitem{fB94}
F.~Borceux, \emph{Handbook of categorical algebra 1}, Encyclopedia of
  Mathematics and Its Applications, vol.~50, Cambridge University Press, 1994.

\bibitem{CM95}
G.~Carlsson and R.J. Milgram, \emph{Stable homotopy and iterated loop spaces},
  in: Handbook of Algebraic Topology (Ed.~I.M. James, ed.), North-Holland,
  1995, pp.~505--583.

\bibitem{CM57}
H.S.M. Coxeter and W.O.J. Moser, \emph{Generators and relations for discrete
  groups}, Ergebnisse der Mathematik und ihrer grenzgebiete, vol.~14,
  Springer-Verlag, 1957.

\bibitem{CY981}
L.~Crane and D.~Yetter, \emph{Deformations of (bi)tensor categories}, Cahier de
  Topologie et G\'eometrie Differentielle Cat\'egorique \textbf{39} (1998),
  163--180.

\bibitem{vD90}
V.~Drinfeld, \emph{Quasihopf algebras}, Leningrad J. Math. \textbf{1} (1990),
  1419--1457.

\bibitem{jE1}
J.~Elgueta, \emph{Cohomology and deformation theory of monoidal 2-categories
  i}, to appear in Adv. Math. (arXiv: math.QA/0204099).

\bibitem{mG64}
M.~Gerstenhaber, \emph{On the deformations of rings and algebras}, Ann. of
  Math. \textbf{79} (1964), 59--103.

\bibitem{sL02}
S.~Lack, \emph{A quillen model structure for 2-categories}, K-Theory
  \textbf{26} (2002), 171--205.

\bibitem{sM63}
S.~MacLane, \emph{Natural associativity and commutativity}, Rice Univ. Studies
  \textbf{49} (1963), 28--46.

\bibitem{sM98}
\bysame, \emph{Categories for the working mathematician}, Third Edition, GTM,
  vol.~5, Springer, 1998.

\bibitem{MS94}
M.~Markl and J.~Stasheff, \emph{Deformation theory via deviations}, J. Algebra
  \textbf{170} (1994), 122--155.

\bibitem{rjM66}
R.J. Milgram, \emph{Iterated loop spaces}, Ann. of Math. \textbf{84} (1966),
  386--403.

\bibitem{rS95}
R.~Street, \emph{Descent theory}, notes of lectures presented at Oberwolfach,
  September 1995 (http://www.maths.mq.edu.au/~street/Descent.pdf).

\bibitem{rS87}
\bysame, \emph{The algebra of oriented simplexes}, J. Pure Appl. Algebra
  \textbf{49} (1987), 283--335.

\bibitem{rS03}
\bysame, \emph{Categorical and combinatorial aspects of descent theory},
  preprint (arXiv: math.CT/0303175) (2003).

\bibitem{cW94}
C.~Weibel, \emph{An introduction to homological algebra}, Cambridge studies in
  advanced mathematics, vol.~28, Cambridge University Press, 1994.

\bibitem{dY98}
D.~Yetter, \emph{Braided deformations of monoidal categories and vassiliev
  invariants}, in: Higher Category Theory (M.~Kapranov E.~Getzler, ed.),
  American Mathematical Society, 1998, A.M.S. Contemporary Mathematics, vol.
  230, pp.~117--134.

\bibitem{dY01}
\bysame, \emph{Functorial knot theory}, Series on Knots and Everything,
  vol.~26, World Scientific, 2001.

\end{thebibliography}
\end{document}